\documentclass[a4paper,11pt]{article}
\usepackage[english]{babel}
\usepackage[latin1]{inputenc}
\usepackage[T1]{fontenc}
\usepackage{amsmath}
\usepackage{amsfonts}
\usepackage{dsfont}
\usepackage{amsthm}
\usepackage{amssymb}
\title{Hyperbolic measure of maximal entropy for generic rational maps of $\P^k$}
\usepackage{dsfont}
\author{Gabriel Vigny}
\begin{document}
\newtheorem{Theorem}{Theorem}
\newtheorem*{Theorem*}{Theorem}
\newtheorem{theorem}{Theorem}[section]
\newtheorem{proposition}[theorem]{Proposition}
\newtheorem{defi}[theorem]{Definition}
\newtheorem{corollaire}[theorem]{Corollary}
\newtheorem{Hypothesis}[theorem]{Hypothesis}
\newtheorem{lemme}[theorem]{Lemma}
\newtheorem{Remark}[theorem]{Remark}
\newcommand{\U}{\mathcal{U}}
\newcommand{\C}{\mathcal{C}}
\renewcommand{\P}{\mathbb{P}}
\renewcommand{\phi}{\varphi}
\newcommand{\Cc}{\mathbb{C}}
\newcommand{\Nn}{\mathbb{N}}
\newcommand{\Rr}{\mathbb{R}}
\newcommand{\Qq}{\mathbb{Q}}
\newcommand{\Zz}{\mathbb{Z}}
\newcommand{\Acal}{\mathcal{A}}
\newcommand{\Bcal}{\mathcal{B}}
\newcommand{\Dcal}{\mathcal{D}}
\newcommand{\Ecal}{\mathcal{E}}
\newcommand{\Gcal}{\mathcal{G}}
\newcommand{\Hcal}{\mathcal{H}}
\newcommand{\Lcal}{\mathcal{L}}
\newcommand{\Ical}{\mathcal{I}}
\newcommand{\Pcal}{\mathcal{P}}
\newcommand{\Qcal}{\mathcal{Q}}
\newcommand{\Scal}{\mathcal{S}}
\newcommand{\Zcal}{\mathcal{Z}}
\hyphenation{plu-ri-sub-har-mo-nic}
\date{}

\maketitle
\begin{abstract}
Let $f$ be a dominant rational map of $\P^k$ such that $\exists s <k, \ \mathrm{with} \ \lambda_s(f)>\lambda_l(f)$ for all $l$. Under mild hypotheses, we show that, for $A$ outside a pluripolar set of $\mathrm{Aut}(\P^k)$, the map $f\circ A$ admits a hyperbolic measure of maximal entropy $\log \lambda_s(f)$ with explicit bounds on the Lyapunov exponents. In particular, the result is true for polynomial maps hence for the homogeneous extension of $f$ to $\P^{k+1}$.
This provides many examples where non uniform hyperbolic dynamics is established.

One of the key tools is to approximate the graph of a meromorphic function by a smooth positive closed current. This allows us to do all the computations in a smooth setting, using super-potentials theory to pass to the limit. 
 \end{abstract}


\section{Introduction}
Let $f:\P^k\to \P^k$ be a dominant meromorphic map of the projective space $\P^k$ (i.e. a rational map). We are interested in the ergodic properties of $f$. More precisely, we want to construct a measure of maximal entropy that we want to compute and then study its properties (ergodicity, mixing, hyperbolicity ... ). This is a natural yet difficult question in dynamics and the tools of complex analysis and complex geometry often allow to answer that question more easily. \\

Such study starts with the computation of the \emph{dynamical degrees}.  For $0\leq l \leq k$, let $L_l$ and $L_{k-l}$ be generic linear subspaces of $\P^k$ of codimension $l$ and $k-l$. Then the number:
$$\lambda_l(f):=\mathrm{Card}(f^{-1}(L_l) \cap L_{k-l})=\mathrm{Card}(L_l \cap f( L_{k-l})) $$
is well defined and does not depend on the choice of $L_i$ as it is defined in cohomology. In particular, if $\omega$ denotes the Fubini-Study form on $\P^k$ then we also have:
$$\lambda_l(f):=\int_{\P^k} f^{*}(\omega^l) \wedge \omega^{k-l}=\int_{\P^k} \omega^l \wedge f_*(\omega^{k-l}).$$
The sequence $(\lambda_l(f^n))$ satisfies $\lambda_l(f^{n+m})\leq \lambda_l(f^n)\lambda_l(f^m)$ so we can define
the $l$-th \emph{dynamical degree} as (see \cite{RS}):
$$d_l(f):= \lim_{n \rightarrow + \infty} (\lambda_l(f^n))^{1/n}.$$
The degree $d_l$ measures the asymptotic spectral 
radius of the action of $f^*$ on the cohomology group $H^{l,l}(\P^k)$. When $\lambda_l(f^n)=\lambda_l(f)^n$ for all $n$, we say that $f$ is $l$-\emph{algebraically stable} (\cite{DS6}).   
The last degree $d_k$ is the \emph{topological degree}. The sequence of degrees is increasing up to a rank $s$ and then it is decreasing (see \cite{Gromov}). \\

Assume that one of the dynamical degree $d_s$ of $f$ is greater than the others. Such map is said to be \emph{cohomologically hyperbolic}. It is conjectured  (see \cite{Gu1}) that there exists a measure of maximal entropy $\log d_s$. That measure should be hyperbolic (no Lyapunov exponent is zero) and the saddle points should be equidistributed along that measure (that last point is out of the scope of the article). Such statement has been proved in the cases where the highest dynamical degree is the topological degree (\cite{Gu2}), for Hénon mappings (see e.g \cite{BS}), regular birational maps (\cite{DS10}), polynomial-like and horizontal-like maps (\cite{DS3, DinhNguyenSibony1}) ... That gives large families of examples where hyperbolic dynamics is established. Still, the result is not known is general and there are natural families for which it is left to be done: birational mappings of $\P^2$, polynomial mappings of $\Cc^k$ ($k\geq 3$), and more generally rational mappings of $\P^k$ for which the highest dynamical degree is not the topological degree (and from now on, that will be the case we are in).  \\

A fruitful approach toward that direction has been initiated by Bedford and Diller for a birational map $f$ on a projective surface $X$ in \cite{BD1}. They define a geometric condition on the indeterminacy sets $I(f)$ and $I(f^{-1})$ under which they can construct the wanted measure (the computation of the entropy was done in \cite{Duj2}). Then they show that, in the case where $X$ is $\P^2$ that condition is generic in the following sense: for any $f$ satisfying that condition and any $A$ outside a pluripolar set of $\mathrm{Aut}(\P^2)$ then $f \circ A$ also satisfies that geometric condition. In \cite{Dilgu}, the authors showed that there exist examples that do not satisfy that condition and gave a more general condition that is still not always satisfied.  Finally, in the articles \cite{DDG2, DDG3}, the authors generalize that idea to the case of a meromorphic map of a projective surface (under a more general integral condition); whereas that gives new families where the program is fulfilled (notably polynomial mappings of $\Cc^2$) it is yet not general. Indeed, a recent work of Buff gives examples where that condition is not satisfied (\cite{Buff}). Getting more and more general conditions in hopping to finally get all the existing meromorphic maps seems to be a failing approach as one always seems to find maps that are a "little bit more pathologic" (that might simply be due to the fact that the above conjecture is false). Still, that approach gives large families of map for which we understand fairly well the chaotic dynamics. Furthermore, pluripolar sets are of zero Lebesgue measure, so for a given map $f$, though we may not be able to construct the right measure for $f$, we are able to do so for arbitrarily small approximations of $f$ (that is a map $f\circ A_\varepsilon$ where $A_\varepsilon$ is close to the identity in $\mathrm{Aut}(\P^k)$).\\

This was one of the motivations of De Th\'elin and the author in \cite{DV1} where we were interested in dynamics in higher dimension. We considered the family of birational maps $f$ of $\P^k$ such that $\mathrm{dim}(I(f))=k-s-1$ and $\mathrm{dim}(I(f^{-1}))=s-1$, for some $1 \leq s\leq k-1$ (when $k=2$, that gives every birational maps but the situation is more complex when $k\geq 3$). We gave a geometric condition on $I(f)$ and $I(f^{-1})$ analogous to Bedford-Diller's condition under which we constructed a measure of maximal entropy. Then, we showed that for any $A$ outside a pluripolar set of $\mathrm{Aut}(\P^k)$, $f \circ A$ satisfies that condition (we do not need that $f$ itself satisfies the condition). A natural question is to prove the same statement for rational maps (not necessarily birational) with no hypothesis on the dimension of the indeterminacy sets. \\

This is exactly the aim of the article. A difference is that we no longer look for a condition that ensures the existence of the right measure, we directly try to construct the measure and we show we can succeed outside a pluripolar set. We denote by $\C_q$ the convex cone of positive closed currents of bidegree $(q,q)$ and mass $1$. The main results of the article can be summed up in the following theorem (see below for notions related to super-potentials theory):
\begin{Theorem}\label{principal} Let $f$ be a dominant rational map of $\P^k$.
\begin{enumerate}
\item Outside a countable union of analytic sets of $A\in\mathrm{Aut}(\P^k)$, the map $f_A:= f\circ A$ satisfies $\lambda_s(f_A^n)=\lambda_s(f)^n$ for all $n$ and $s$.
\item Assume that $\exists s <k \ \mathrm{with} \ \lambda_s(f)>\lambda_l(f)$ for all $l<s$. 
Then outside a pluripolar set of $A\in\mathrm{Aut}(\P^k)$,
for any smooth form $\Omega_s \in \C_s$, the sequence of currents $ \lambda_s(f)^{-n}(f^n_A)^*(\Omega_s)$
converges in the Hartogs' sense to the Green current $T^+_{s,A}$ which is $f_A^*$-invariant.
\item Assume that $\exists s <k, \ \mathrm{with} \ \lambda_s(f)>\lambda_l(f)$ for all $l$. 
Then outside a pluripolar set of $A\in\mathrm{Aut}(\P^k)$,
for any smooth form $\Omega_{k-s} \in \C_{k-s}$, the sequence of currents $\lambda_s(f)^{-n} (f^n_A)_*(\Omega_{k-s})$
converges in the Hartogs' sense to the Green current $T^-_{s,A}$ which is $(f_A)_*$-invariant. Furthermore, the measure $\nu_A:=T^+_{s,A} \wedge T^-_{s,A}$ is well defined in the sense of super-potentials.
\item If in addition the map $f$ satisfies $\mathrm{dim}(I(f))=k-s-1$  or $I(f)\subset H$ for a hyperplane $H$ then outside a pluripolar set of $\P^k$ the measure  $\nu_A$ is an invariant measure of maximal entropy $\log \lambda_s(f)$ which is hyperbolic.
\item Assume that $f$ is a polynomial map of $\Cc^k$, then the points 1, 2, 3 and 4 are true replacing $\mathrm{Aut}(\P^k)$ with $\mathrm{Aff}(\Cc^k)$, the affine automorphisms of $\Cc^k$.  
\end{enumerate}
\end{Theorem}
An important remark is that for polynomial maps, the inderminacy set is always contained in the hyperplane
at infinity so point 4 and 5 holds for polynomial mappings as soon as  $\exists s <k, \ \mathrm{with} \ \lambda_s(f)>\lambda_l(f)$. Then starting with a rational map $f$ of $\P^k$, though we might not have point 4 in Theorem \ref{principal} for $f$, we do have it for any homogeneous extension $\widetilde{f}$ of $f$ to $\P^{k+1}$. Since $f$ is a factor of $\widetilde{f}$, it means that though we might not be able to approximate $f$ by hyperbolic maps in the orbit under $\mathrm{Aut}(\P^k)$, we can approximate a more complex dynamics ($\widetilde{f}$) but in a bigger space ($\mathrm{Aut}(\P^{k+1})$ or  $\mathrm{Aff}(\Cc^{k+1})$). 

Observe also that the case of birational maps of $\P^{3}$ is covered by Theorem \ref{principal}: as $\mathrm{dim}(I(f^{\pm 1})) \leq 1$, if $\lambda_1(f)>\lambda_2(f)$ then we apply it directly,  if $\lambda_2(f)>\lambda_1(f)$ then we apply it to $f^{-1}$.    \\

Let us explain the approach of the article and how it differs from the one in \cite{DV1}.
For birational maps of \cite{DV1}, the hypothesis on the dimension of the indeterminacy sets
implies that if $\cup_{i\in \Nn} (f^{-1})^i(I^+) \cap I^-=\varnothing$ then the maps is algebraically stable, $d_s=d^s=\delta^{k-s}$ is the highest dynamical degree and $f^*(\omega^s)=f^*(\omega)^s$ ($d$ is the algebraic degree of $f$, $\delta$ the algebraic degree of $f^{-1}$). We then considered the following condition on $f$:
\begin{equation}\label{condition_analytique}
\begin{cases} 
\sum_{n\in\Nn} \frac{1}{d^{sn}}\int_{f^{n}(I^-)} u f^*(\omega)^{s-1} >-\infty\\
\sum_{n\in\Nn} \frac{1}{\delta^{(k-s)n}}\int_{f^{-n}(I^+)} u' f_*(\omega)^{k-s-1} >-\infty,
\end{cases}
\end{equation}
where $u$ (resp. $u'$) is a quasi-potential of $f^*(\omega)$ (resp. $f_*(\omega)$).
Then, under (\ref{condition_analytique})
we were able to construct the Green currents and a (mixing) hyperbolic measure of maximal entropy using Theorem 1 in \cite{DV1} and the results of \cite{DT1}. 

We then showed that (\ref{condition_analytique}) is given by a decreasing sequence of some quasi-plurisubharmonic function $g_n$ on $\mathrm{Aut}(\P^k)$ and then we provide examples in the orbit of $f$ to show that $g:= \lim_n g_n \not\equiv -\infty$ (these examples also showed that the condition giving algebraic stability was satisfied outside a countable union of analytic sets). To give the expression of such function $g_n$, we see it as the push-forward on $\P^k$ of some current in $\mathrm{Aut}(\P^k)\times \P^k$. A key point to make the computations was that, thanks to the hypothesis on the dimension of $I^+$ and $I^-$, such current was smooth
outside a set of codimension $2$ in $\mathrm{Aut}(\P^k)$ (Lemma 3.3.2) hence the computation of the $dd^c$ of the currents was just the trivial extension of its $dd^c$ wherever it is well defined (\cite{dem}[Chapter III Corollary 4.11]). 

Dealing with higher dimension with no control on $I(f)$ is the main difficulty of the article. In order to deal with such currents, we use the theory of super-potentials of Dinh and Sibony (\cite{DS6}). Though the general strategy is similar, some serious obstructions appear that force us to make deep changes. First of all, the indeterminacy sets $I$ and $I'$ (see the definition later on though at this point the reader may think of them as $I^+$ and $I^-$) do intersect in general. Hence, algebraic stability cannot be given by a simple condition of the 
form of "$\cup_{i\in \Nn} (f^{-1})^i(I^+) \cap I^-=\varnothing$" (they should be such  condition taking into account that the intersection of the different stratifications of $I^-$ and their images with $I^+$ are transverse but it would not of any use). Instead, we proceed inductively and show that, under algebraic conditions, we can defined $(f^n)^{*}([M])$ with $(f^n)^{*}([M])=(f^*)^{n}([M])$ for a given analytic set $M$ of codimension $s$.
Providing explicit examples where algebraic stability holds (the spirit of such examples follows ideas of Dinh), we then show that outside a countable union of analytic sets of $A\in\mathrm{Aut}(\P^k)$, pull-backs and push-forwards are well defined in the sense of super-potentials.

 In a second part, we construct the Green currents and their intersection. 
 Instead of giving a condition (\ref{condition_analytique}), we directly try to construct the measures and currents, we show we can succeed outside a pluripolar set.   For that, the idea is to consider the rational map:
\begin{align*}
\widetilde{F}: \mathrm{Aut}(\P^k)\times \P^k &    \to \mathrm{Aut}(\P^k) \times \P^k\\
(A,z) &\mapsto  (A,f_A(z)),
\end{align*}
and to show that $d_s(f)^{-n}(\widetilde{F}^*)^n(\omega^s)$ is well defined and that its slices converges (outside a pluripolar set of $\mathrm{Aut}(\P^k)$) to the Green current of $f_A$ in the sense of super-potentials. For that we want to compute the value of the slice of a quasi-potential of $d_s(f)^{-n}(\widetilde{F}^*)^n(\omega^s)$ at a smooth form of $\C_{k-s+1}$. Then we want to show that it defines a DSH function computing its $dd^c$ and providing examples where it is finite (using the same kind of examples as above).  The difficulty lies in the fact that such $dd^c$ is not a priori clearly defined since we have no control on the singularities of $\widetilde{F}$ (contrary to \cite{DV1}). To overcome that problem, we regularize the map $\widetilde{F}$ in the following sense: we approximate its graph by a smooth positive closed current. Though we do not have a map anymore, we preserve the cohomology and we keep the functional properties of the pull-back and push-forward. Then all the computations make sense and we pass to the limit for $\widetilde{F}$ using pluripotential theory. We believe that idea can be used in other cases.   
In a last section, we prove points 4 and 5 in Theorem \ref{principal}. We use Theorem 1 in \cite{DV1} to show that the entropy of $f_A$ is $\log \lambda_s(f)$, the hyperbolicity is obtained thanks to the results of \cite{DT1}. As above, the idea is to prove that the desired properties are obtained under DSH conditions. We need the additional hypotheses of point 4 on the indeterminacy sets to construct examples that satisfy these conditions. In an independent paragraph, we explain how knowing the entropy and hyperbolicity of the homogeneous extension $\widetilde{f}$ gives the entropy and hyperbolicity of $f$ using the theory of the entropy of a skew-product. \\

\noindent {\bf Acknowledgements.} I am grateful to De Thélin for numerous conversations where he convinced me that the results of the paper were achievable and for explaining how Corollary 3 in \cite{DT1} could be used here. \\

\noindent {\bf Notations and preliminaries.}
In what follows, $f:\P^k\rightarrow\P^k$ denotes a meromorphic map. Such a map is holomorphic
outside an analytic subset $I(f)$ of codimension $\geq 2$ in $\P^k$.  It can be written in homogeneous coordinates as $[P_0: \dots:P_k]$ where the $P_i$ are homogeneous polynomials of algebraic degree $d$ in the $(z_0,\dots,z_k)$ variable, with $\gcd_i(P_i)=1$.
Let $\Gamma$ denote the
closure of the graph of the restriction of $f$ to $\P^k\setminus
I(f)$. This is an irreducible analytic set of dimension $k$ in
$\P^k\times\P^k$. Let $\pi_1$ and $\pi_2$ denote the canonical
projections of $\P^k\times \P^k$ on the factors. 
The indeterminacy set  $I(f)$ is also the set of points $z\in\P^k$ such that $\dim
\pi_1^{-1}(z)\cap\Gamma\geq 1$. We sometimes write $I$ instead of $I(f)$. We assume that $f$ is {\it dominant},
that is, $\pi_2(\Gamma)=\P^k$. The {\it second indeterminacy set} of
$f$ is the set $I'$  of points $z\in\P^k$ such that $\dim
\pi_2^{-1}(z)\cap\Gamma\geq 1$. Its codimension is also at least equal
to $2$.
If $A$ is a subset of $\P^k$, define
$$f(A):=\pi_2(\pi_1^{-1}(A)\cap\Gamma)\quad \mathrm{and}\quad
f^{-1}(A):=\pi_1(\pi_2^{-1}(A)\cap\Gamma).$$
We will need to distinguish between the direct image of $A$ by $f$ iterated $n$ times (that we denote $(f^n)(A)$) and the direct image iterated $n$ times of $A$ by $f$ (that we denote $(f)^n(A)$). We use the same notations for preimages. If $f$ is holomorphic, both notions coincide. That does not need to be the case if $f$ is meromorphic. \\

We need to define pull-back and push-forward of positive closed currents. Recall that if $S$ is a positive closed current of bidegree $(s,s)$, we denote its \emph{mass} $\|S\|:= \int S \wedge \omega^{k-s}$.
Define formally for a current $S$ on $\P^k$, not necessarily positive
or closed, the pull-back $f^*(S)$ by
\begin{equation} \label{eq_pullback_def}
f^*(S):=(\pi_1)_*\big(\pi_2^*(S)\wedge [\Gamma]\big).
\end{equation}
This 
makes sense if the wedge-product $\pi_2^*(S)\wedge
[\Gamma]$ is well defined, in particular, when $S$ is smooth. We will be particularly interested in the case 
where $S$ is the current of integration on ananalytic set.
Similarly, the operator $f_*$ is formally defined by 
\begin{equation}  \label{eq_pushforward_def}
f_*(R):=(\pi_2)_*\big(\pi_1^*(R)\wedge [\Gamma]\big).
\end{equation}
We need in the article the theory of super-potentials (\cite{DS6} for proofs, or the appendix of \cite{DV1}). The formalism of super-potentials allows to extend the calculus of potentials to the case of general bidegree. Recall that if $T\in \C_q$, it is cohomologous to a fixed smooth form $\Omega_q\in \C_q$, hence we can write it $T= \Omega_q+dd^c U_T$ where $U_T$ is a quasi-potential of $T$. A super-potential $\U_T$ of $T$ is then the function defined for smooth $ S\in \C_{k-q+1}$ by $\U_T(S)=\langle U_T, S\rangle$. This definition can be extended to arbitrarily elements of $\C_{k-q+1}$ making $\U_T$ a quasi-plurisubharmonic function on $\C_{k-q+1}$ (according to the notion of \emph{structural variety} on $\C_{k-q+1}$).  

In particular, the notion of pull-back and push-forward can be extended to $f^*$-admissible elements of $\C_q$ (resp. $f_*$-admissible elements of $\C_{k-q}$) that is elements whose super-potentials are finite at  $\lambda(f_{s-1})f_*(R)$ for some $R$ smooth in $\C_{k-s+1}$ (resp. at  $\lambda(f_{s+1})f_*(R)$ for some $R$ smooth in $\C_{s+1}$). For smooth forms in $\C_q$, the notions of pull-back and push-forward coincide with the ones given by (\ref{eq_pullback_def}) and (\ref{eq_pushforward_def}) and super-potentials extend that notion to admissible elements by pluri-subharmonicity along the structural varieties.

Finally, recall that the group $\mathrm{Aut}(\P^k)$ of automorphisms of $\P^k$ is $\mathrm{PGl}(\Cc^{k+1})$. In particular, it is a Zariski dense open set in $\P^{(k+1)^2-1}(\Cc)$. For $A \in \mathrm{Aut}(\P^k)$ and $l\leq k$, one has that $\lambda_l(f\circ A)= \lambda_l(A \circ f)=\lambda_l(f)$ (the quantity $\lambda_l(f)$ is defined in the beginning of the introduction). 
  This explain why we choose to consider such perturbations of $f$. Indeed, it could seem that a natural way to approximate $f$ would be to slighty change the polynomials $P_i$ (where $f=[P_0:\dots:P_k]$). But such perturbation gives generically a holomorphic map as the common zero set of $k+1$ polynomials in $\P^k$ is generically empty. On the contrary, our choice ensures that we stay in the same family. 
  
    For $h$ in the orbit of $f$, we use the notation $L_h:=\lambda_q(f)^{-1} h^*$ (resp. $\Lambda_h:=(\lambda_{k-q})^{-1}h_*=(\lambda^-_q(f))^{-1} h_*$) for the normalized pull-back (resp. push-forward) in the sense of super-potentials acting on  $h^*$-admissible elements of $\C_q$ (resp. $h_*$-admissible elements of $\C_{k-q}$). We simply write $L$ and $\Lambda$ instead of $L_f$ and $\Lambda_f$. \\

\section{Algebraic stability is dense}\label{Algebraic stability is dense}
The purpose of the section is to prove the following theorem (which gives the first point of Theorem \ref{principal}). We do not assume here that $\lambda_s(f)$ is greater than the other $\lambda_l(f)$. The results of the section remain true for a meromorphic correspondence but we state them in the case of a meromorphic map for simplicity. Let $\Omega_q\in \C_q$ be a fixed smooth form, we denote by $\U_{L_h(\Omega_q)}$ (resp. $\U_{\Lambda_h(\Omega_q)}$) a super-potential of $L_h(\Omega_q)$ (resp. $\Lambda_h(\Omega_q)$).
\begin{theorem}\label{algebraic_stability_generic}
For all $n$, there exists a Zariski dense open set $\mathcal{Z}_{n,s}$ of elements $h$ in the orbit of $f$ for which:
\begin{enumerate}
 \item $\lambda_s(h^m)=\lambda_s(h)^m=\lambda_s(f)^n$ for all $m\leq n$;
 \item $(h^*)^m=(h^m)^*$ and $(h_*)^m=(h^m)_*$ for all smooth forms in $\C_s$ and $\C_{k-s}$ in the sense of super-potentials;  
 \item If $\U_T$ is a super-potential of $T$ smooth in $\C_s$, then the following $\U_{L_h^n(T)}$ is a super-potential of $L_h^n(T)$ on smooth forms:
\begin{align}\label{super-potential_L_h^n}
\U_{L_h^n(T)}&= \sum_{i\leq n} \left( \frac{d_{s-1}}{d_s}\right)^i \U_{ L_h(\Omega_s)}\circ \Lambda_h^{i} +   \left(\frac{d_{s-1}}{d_s}\right)^n \U_T\circ \Lambda_h^{n}; 
\end{align}
\item If $\U_S$ is a super-potential of $S$ smooth in  $\C_{k-s}$, then the following $\U_{\Lambda_h^n(S)}$ is a super-potential of $\Lambda_h^n(S)$ on smooth forms:
\begin{align}\label{super-potential_Lambda_h^n}
\U_{\Lambda_h^n(S)}&= \sum_{i\leq n} \left( \frac{d_{s+1}}{d_s}\right)^i \U_{ \Lambda_h(\Omega_s)}\circ L_h^{i} +   \left(\frac{d_{s+1}}{d_s}\right)^n \U_S\circ L_h^{n}. 
\end{align}
 \end{enumerate}
Furthermore, the intersection $\cap_{n \in \Nn}\mathcal{Z}_{n,s}$ contains an open set $\mathcal{Y}$.
\end{theorem}
We denote $\mathrm{dim}(I):=m$ and $\mathrm{dim}(I'):=m'$. We consider the set $\C_1:= \pi_2^{-1}(I')\cap \Gamma$ (it is the critical set for $(\pi_2)|_{\Gamma}$, the second projection of $(\P^{k})^2$, and an exceptional set of $\Gamma$). It is an analytic
subset of $\Gamma$ so it has dimension $\mathrm{dim}(\C_1) \leq k-1$. Similarly, we consider $\C'_1:= \pi_1^{-1}(I)\cap \Gamma$ which is an analytic set of dimension $\leq k-1$. In particular, for $p \in I'$ generic, we have $\mathrm{dim}(\pi_2^{-1}(p)\cap \Gamma) =\mathrm{dim}(\C_1)-m' \leq k-1-m'$. 
For $r \in \{ \mathrm{dim}(\C_1)-m', \dots, \mathrm{dim}(\C_1)\}$, we consider:
 $$I'_r:= \{p \in I', \ \mathrm{dim}(\pi_2^{-1}(p)\cap \Gamma) =r  \}.$$    
Then $I'_r$ is a (possibly empty) locally analytic set of dimension $\leq  \mathrm{dim}(\C_1)-r$ (which is less than $k-1-r$)
and $ \cup_{r'\geq r} I'_{r'}$ is an analytic set. Similarly, for $r \in \{ \mathrm{dim}(\C'_1)-m, \dots, \mathrm{dim}(\C'_1)\}$, we consider:
 $$I_r:= \{p \in I, \ \mathrm{dim}(\pi_1^{-1}(p)\cap \Gamma) =r  \}.$$    
Then $I_r$ is a (possibly empty) locally analytic set of dimension $\leq  \mathrm{dim}(\C'_1)-r$ (which is less than $k-1-r$)
and $ \cup_{r'\geq r} I_{r'}$ is an analytic set. 

Recall that if $M$ is an analytic set then for any $\varepsilon >0$, there exists a $\delta >0$ such that if $U_\delta$ is a $\delta$-neighborhood of $M$ then $f^{-1}(U_\delta)$ is contained in a $\varepsilon$-neighborhood of $f^{-1}(M)$. The same result holds for direct image.  
\begin{lemme}\label{pull_back_generic_linear}
Let $M$ be an analytic set of codimension $s$ such that for all $r \in  \{ \mathrm{dim}(\C_1)-m', \dots, \mathrm{dim}(\C_1)\}$, $M \cap \cup_{r'\geq r} I'_{r'}$ is empty if $ \dim (\cup_{r'\geq r} I'_{r'}) \leq s-1$ and of dimension $\dim (\cup_{r'\geq r} I'_{r'}) - s$ if not.
 Assume also that no component of $\pi_2^{-1}(M)\cap \Gamma$ is contained in $\C'_1$. 
\begin{enumerate}
\item Then $f^{-1}(M)$ is an analytic set of codimension $s$ such that $\mathrm{codim} (f^{-1}(M) \cap I(f)) \geq s+1$. For all analytic set $M' \subset M$ of codimension $\geq s+1$, then $\mathrm{codim}(f^{-1}(M'))\geq s+1$.  
\item The positive closed current $f^*([M])$ of bidegree $(s,s)$ is well defined, depends continuously on $[M]$ in the sense of currents and is equal to $ [f^{-1}(M)]$ (counting with multiplicity). Hence, we have that $\lambda_s(f) = \| f^*([M]) \| \times \|[M]\| ^{-1}$.  
\end{enumerate}
\end{lemme}
\noindent \emph{Proof.} Take $M$ as in the lemma, we prove the first point. We have that $\pi^{-1}_2(M)  \cap \Gamma$
 is an analytic set which is of codimension $s$ outside $\C_1$. For $r \geq \mathrm{dim}(\C_1)-s+1$, we have that $\mathrm{dim}(I'_r)<s $ hence $ \pi_2^{-1}(M \cap I'_{r})=\varnothing$. Now, for $r$ such that $ \mathrm{dim}(\C_1)-m' \leq r \leq  \mathrm{dim}(\C_1)-s$, we have that:
$$ \mathrm{dim}( \pi_2^{-1}(M \cap I'_{r}) \cap \Gamma)  \leq  \mathrm{dim}(\C_1)-r+k-s-k+r\leq k-s-1 .$$ 
That implies that $\pi_2(M)  \cap \Gamma$ is an analytic set of codimension $s$. Pushing-forward by $\pi_1$ (and keeping track of the multiplicity), we have that $f^{-1}(M)$ is indeed an analytic set of dimension $s$ as $\pi_2(M)  \cap \Gamma \nsubseteq \C'_1$ and $\mathrm{codim} (f^{-1}(M) \cap I(f)) \geq s+1$. 
For the second part of that point, take $M'$ of codimension $\geq s+1$. Then, outside $I(f)$, it is clear that $\mathrm{codim} (f^{-1}(M')) \geq s+1$ and the previous argument shows that $\mathrm{codim}(f^{-1}(M') \cap I(f)) \geq  \mathrm{codim}(f^{-1}(M) \cap I(f)) \geq s+1$.\\ 

Now, for the second point, we have that $[\pi_2^{-1}(M)\cap \Gamma]$ is a well defined positive closed current of 
bidegree $(s,s)$ as the current of integration on an analytic set of codimension $s$. Consider $\Gamma' := \Gamma \setminus \left( \C'_1 \cup \C_1 \right)$, it is a complex manifold as the graph of a map.  The first point of the lemma gives that $[\pi_2^{-1}(M)\cap \Gamma]$ is equal to the trivial extension of $[\pi_2^{-1}(M)\cap \Gamma']$ (because $\pi_2^{-1}(M)\cap (\C_1 \cup \C'_1)$ is of codimension $\geq s+1$ and both currents coincide outside a set of zero mass for them.). Furthermore, the fiber of $\pi_2$ restricted to $\Gamma'$ are either finite sets or empty. Theorem 1.1 in \cite{DS13} implies that $(\pi_2)_{|\Gamma'}^*([M])=[\pi_2^{-1}(M)\cap \Gamma']$ is a well defined positive closed current on $\Gamma'$ (that is it depends continuously on $[M]$ for the topology of current: if $(R_n)$ is a sequence of smooth currents converging to $[M]$ then $(\pi_2)_{|\Gamma'}^*(R_n)$ converges to $(\pi_2)_{|\Gamma'}^*([M])$). In particular, $[\pi_2^{-1}(M)\cap \Gamma]$ depends continuously on $M$. Pushing-forward by $(\pi_1)_*$ gives that $f^*([M])$ is well defined, depends continuously on $M$ in the sense of currents and is equal to $[f^{-1}(M)]$ (again, we keep track of the multiplicity). 
Now we deduce that  $\lambda_s(f) = \| f^*([M]) \| \times \|[M]\| ^{-1}$ (that would be true if $[M]$ was a smooth current and we can conclude by continuity). \hfill $\Box$ \hfill   \\

Similarly, one can prove:
\begin{lemme}\label{push_forward_generic_linear}
Let $N$ be an analytic set of dimension $s$ such that for all $r \in  \{ \mathrm{dim}(\C'_1)-m, \dots, \mathrm{dim}(\C'_1)\}$, $N	 \cap \cup_{r'\geq r} I_{r'}$ is empty if $ \dim (\cup_{r'\geq r} I_{r'}) \leq k-s-1$ and of dimension $s+\dim (\cup_{r'\geq r} I_{r'}) - k$ if not. Assume also that no component of $\pi_1^{-1}(N)\cap \Gamma$ is contained in $\C_1$.
\begin{enumerate}
\item Then $f(N)$ is an analytic set of dimension $s$ such that $\dim (f(N) \cap I'(f)) \leq s-1$. For all analytic set $N' \subset N$ of dimension $\leq s-1$, then $\mathrm{dim}(f(N'))\leq s-1$.  
\item The positive closed current of bidegree $(k-s,k-s)$ $f_*([N])$ is well defined, depends continuously on $[N]$ in the sense of currents and is equal to $ [f(N)]$ (counting the multiplicity). Furthermore, we have that $\lambda_{s}(f) = \| f_*([N]) \| \times \|[N]\| ^{-1}$.  
\end{enumerate}
\end{lemme}

In order to simplify the exposition, we need the following \emph{ad hoc} definition:
\begin{defi}
An analytic set of codimension $s$ (resp. dimension $s$) is said to be $f^*$-compatible (resp. $f_*$-compatible) if it satisfies the hypotheses of Lemma \ref{pull_back_generic_linear} (resp. Lemma \ref{push_forward_generic_linear}).
\end{defi}
A crucial point for a process in dynamics is that it needs to be iterated. Recall that $\Omega_s\in \C_s$ and $\Omega_{k-s+1}\in \C_{k-s+1}$ are fixed smooth elements.
\begin{proposition}\label{iterated_pull_back}
\begin{enumerate}
\item  Let $M$ be an analytic set of codimension $s$. 
Assume that for all $0\leq m\leq n-1$, $(f^{-1})^m(M)$ is $f^*$-compatible. Then, moving $M$ a little, we can assume that it is $(f^n)^*$-compatible and $f^{-n}(M)=(f^{-1})^n(M)$ (counting the multiplicity) up to a set of codimension $\geq s+1$. Consequently, $\lambda_s(f^n) =\lambda_s(f)^n$.
\item The same result holds for direct images, replacing  $f^*$-compatibility with  $f_*$-compatibility.
\item Assume that there also exists  an analytic set $F$ of dimension $s-1$ satisfying $(F\cup f(F)) \cap(\cup_{0\leq m\leq n-1} (f^{-1})^{m}(M))=\varnothing$. Then for any smooth $T \in \C_s$ and $j\leq n-1$, $L^j(T)$ is $f^*$-admissible and $(f^n)^*(T)=(f^*)^n(T)$ in the sense of super-potentials. If $\U_T$ is a super-potential of $T$ smooth, then the following $\U_{L^n(T)}$ is a super-potential of $L^n(T)$ on smooth forms:
\begin{equation*}
\U_{L^n(T)}= \sum_{i\leq n} \left( \frac{d_{s-1}}{d_s}\right)^i \U_{ L(\Omega_s)}\circ \Lambda^{i} +   \left(\frac{d_{s-1}}{d_s}\right)^n \U_T\circ \Lambda^{n}.
\end{equation*}
\item Similarly, let $N$ be an analytic set of dimension $s$ such that 
 for all $0\leq m\leq n-1$, $(f)^m(N)$ is $f_*$-compatible.
Assume that there exists  an analytic set $E$  of codimension $s+1$ satisfying $(E\cup f^{-1}(E)) \cap(\cup_{0\leq m\leq n-1} (f)^{m}(N))=\varnothing$. Then for any smooth $S \in \C_{k-s}$ and $j\leq n-1$, $\Lambda^j(S)$ is $f_*$-admissible and $(f^n)_*(S)=(f_*)^n(S)$ in the sense of super-potentials. If $\U_S$ is a super-potential of $S$ smooth, then the following $\U_{\Lambda^n(S)}$ is a super-potential of $\Lambda^n(S)$ on smooth forms:
\begin{equation*}
\U_{\Lambda^n(S)}= \sum_{i\leq n} \left( \frac{d_{s+1}}{d_s}\right)^i \U_{ \Lambda(\Omega_{k-s})}\circ L^{i} +   \left(\frac{d_{s+1}}{d_s}\right)^n \U_S\circ L^{n}.
\end{equation*}
\end{enumerate}
\end{proposition}
\noindent \emph{Proof.} Lemma \ref{pull_back_generic_linear} gives that $(f^{-1})^n(M)$ is an analytic set of codimension $s$ and mass $ \lambda_s(f)^n  \|[M]\|$.
Since $f^*$-compatibility is generic and depends continuously on analytic sets, we can indeed take an analytic set $M_1$ close to $M$ such that $(f^{-1})^m(M_1)$ is $f^*$-compatible for $m\leq  n-1$ and $(f^n)^*$-compatible. 
Since $\cup_{m\leq n-1} (f)^m(I') \supset I'(f^n)$, we have that 
$$(f^{-1})^{-n}(M_1 \cap (\cup_{m\leq n-1} (f)^m(I'))^c)=(f^n)^{-1}(M_1 \cap (\cup_{m\leq n-1} (f)^m(I'))^c) $$
(with muliplicity). We claim that $M_1 \cap (\cup_{m\leq n-1} (f)^m(I'))$ is of codimension $\geq s+1$. It is by hypothesis for $I'$, we check it for $f(I^-)$. We have that 
$$M_1 \cap f(I') \subset f_{|I^c} (f^{-1}(M_1) \cap I') \cup (M_1 \cap \pi_{2}(\mathcal{C}'_1)).$$   
 Thus, $\mathrm{codim}(M_1 \cap f(I'))\geq k+1$  (the image of a set of codimension $\geq s+1$ by a holomorphic map is again of  codimension $\geq s+1$
and $\mathrm{codim}(M_1 \cap \pi_{2}(\mathcal{C}'_1))\geq s+1$ by hypothesis). The proof is similar for $M_1 \cap  (f)^m(I')$.

Lemma \ref{pull_back_generic_linear} implies that $ (f^{-1})^{n}( M_1 \cap (\cup_{m\leq n-1} (f)^m(I'))$ is of codimension $\geq s+1$. 
Thus $(f^n)^{-1}(M_1)$ coincides with $(f^{-1})^n(M_1)$ outside a set where $(f^{-1})^n(M_1)$ has zero mass. That implies that  $\lambda_s(f^n) \geq \lambda_s(f)^n$. As the other inequality always stands,  the equality $\lambda_s(f^n) =\lambda_s(f)^n$ follows from the last point of Lemma \ref{pull_back_generic_linear} and the first point is proved. The proof of the second point is the same. \\

 We now prove the third point by induction on $n$ (which is clear for $n=0$). So assume the third point is true for $n-1$. We can choose small neighborhoods $U$ of $M$ and $V$ of $F$ such that $(f^{-1})^{m-1}(U) \cap f(V)=\varnothing$. Let $R$ be a smooth element of $\C_{k-s+1}$ with support in $V$ (for example a regularization of $\mathrm{vol(F)}^{-1}[F]$ using an approximation of the identity in $\mathrm{PGL}(\Cc^{k+1})$). 

Any current $T' \in \C_s$ with support in $U$ is such that $(f^{*})^{n-1}(T')$ has support in  $(f^{-1})^m(U)$ hence we can choose a quasi-potential of $(f^{*})^{n-1}(T')$ as a form with $\mathcal{C}^1$ coefficients outside $U$. In particular, its super-potentials  are finite at $\Lambda(R)$.
That gives that the current $L^{n-1}(T')$ is $f^*$-admissible. 

In particular, for $T$ smooth, we have that  $L^{n-1}(T)$ is more $H$-regular than  $L^{n-1}(T')$ hence it is also $f^*$-admissible. If $\U_{n-1}$ denotes a super-potential of 
$L^{n-1}(T)$, then we have that:
$$  \U_{L^n(T)}= \U_{ L(\Omega_s)} +   \frac{d_{s-1}}{d_s} \U_{n-1}\circ \Lambda $$
on smooth forms. A symmetric argument implies that for any smooth form $R\in \C_{k-s+1}$, then $\Lambda^j(R)$ is well defined in the sense of super-potentials for $j\leq n$ (though we do not claim that $(f_*)^j(R)=(f^j)_*(R)$).
In particular, the induction's hypothesis shows that 
$$  \U_{L^n(T)}= \sum_{i\leq n} \left( \frac{d_{s-1}}{d_s}\right)^i \U_{ L(\Omega_s)}\circ \Lambda^{i} +   (\frac{d_{s-1}}{d_s})^n \U_T\circ \Lambda^{n} $$
on smooth forms. The same proof gives the result for direct image.  \hfill $\Box$ \hfill   \\

Taking the intersection over all $n\in \Nn$, the above theorem means that algebraic stability is generic in the orbit of $f$ under $\mathrm{Aut}(\P^k)$ (it stands outside a countable union of analytic varieties). We now provide explicit examples in the orbit of $f$ to show that it is not empty.  \\

Let $E^-_{k-s}$ and $E^+_{s}$ be two linear subspaces of complex dimension $k-s$ and $s$ that are respectively $f^*$-compatible and $f_*$-compatible with $E^-_{k-s} \cap E^+_{s}=\{p\}$ reduced to a point. We can then choose $E^-_{k-s-1}\subset E^-_{k-s}$ and $E^+_{s-1}\subset E^+_{s}$ two linear subspaces of $\P^k$ of complex dimension $k-s$ and $s$ with $p\notin E^+_{s-1}\cup E^-_{k-s-1}$. By Lemmas \ref{pull_back_generic_linear} and \ref{push_forward_generic_linear}, we have that $f^{-1}(E^-_{k-s-1})$ and $f(E^+_{s-1})$ have complex dimension $k-s-1$ and $s-1$.
We claim that we can assume that:
 \begin{align*}
 f^{-1}(E^-_{k-s}) \cap E^+_{s-1}=\varnothing \ \mathrm{and} \ & f(E^+_{s})\cap E^-_{k-s-1}= \varnothing \\
f^{-1}(E^-_{k-s-1})\cap E^+_{s} =\varnothing \ \mathrm{and} \ &f(E^+_{s-1})\cap E^-_{k-s}= \varnothing.
 \end{align*}
Indeed, each of these conditions is generic and can be achieved by moving either $E^+_s$ and $E^+_{s-1}$ or $E^-_{k-s}$ and $E^-_{k-s-1}$. We can choose the homogeneous coordinates $[z_0: \dots: z_k]$ such that
\begin{align*}
E^-_{k-s}=\{z_{0}=\dots= z_{s-1}=0 \}&  \ \mathrm{and} \quad E^+_{s}=\{z_{s+1}=\dots= z_{k}=0 \} \\
E^-_{k-s-1}=\{z_{0}=\dots= z_{s-1}={z_s}=0 \}& \ \mathrm{and} \quad E^+_{s-1}=\{z_{s}=\dots= z_{k}=0 \}.
\end{align*}
In particular, we have:
\begin{itemize}
\item The sets $f^{-1}(E^-_{k-s})$ and $f^{-1}(E^-_{k-s-1})$ (resp. $f(E^+_{s})$ and $f(E^+_{s-1})$ ) are analytic sets of dimension $k-s$ and $k-s-1$ (resp. $s$ and $s-1$). 
\item for every $\varepsilon>0$,
there exist $\delta$-neighborhoods $\mathcal{O}^-$ and $\mathcal{O}^-_1$ (resp. $\mathcal{O}^+$ and $\mathcal{O}^+_1$) of $E^-_{k-s}$ and $E^-_{k-s-1}$
(resp. $E^+_{s}$ and $E^+_{s-1}$) such that $f^{-1}(\mathcal{O}^-)$ and $f^{-1}(\mathcal{O}^-_1)$ (resp. $f(\mathcal{O}^+)$ and  $f(\mathcal{O}^+_1)$) 
are contained in a $\varepsilon$-neighborhood of $f^{-1}(E^-_{k-s})$ and $f^{-1}(E^-_{k-s-1})$ (resp. $f(E^+_{s})$ and $f(E^+_{s-1})$).
 \item finally, choosing $\delta$ small enough,  $\overline{f^{-1}(\mathcal{O}^-)} \cap \overline{\mathcal{O}^+_1} =\varnothing$ and $\overline{f^{-1}(\mathcal{O}^-_1)} \cap \overline{\mathcal{O}^+}=\varnothing$ (resp. $\overline{f(\mathcal{O}^+)} \cap\overline{\mathcal{O}^-_1}=\varnothing$ and $\overline{f(\mathcal{O}^+_1)} \cap \overline{\mathcal{O}^-}=\varnothing$). 
\end{itemize}

Let $A_{\alpha}$ be the element of $\mathrm{Aut}(\P^k)$ defined  by 
\begin{align*} 
&A_{\alpha} ([z_0 : z_1 : \dots z_{s-1} : z_s: z_{s+1} :\dots: z_k ])=  \\
&\quad [ \alpha^{-1} z_0 :\alpha^{-1} z_1 : \dots \alpha^{-1} z_{s-1} : z_s: \alpha z_{s+1} :\dots: \alpha z_k ] 
 \end{align*}     
where $1>\alpha >0$. Choose $\alpha_1$ small enough so that  $A_{\alpha_1}^{-1}(f^{-1}(\mathcal{O}^-)) \Subset \mathcal{O}^-$. This is possible because  $\overline{f^{-1}(\mathcal{O}^-)} \cap \{z_{s}=\dots= z_{k}=0 \} =\varnothing$. Similarly, we can assume that $A_{\alpha_1}^{-1}(f^{-1}(\mathcal{O}^-_1)) \Subset \mathcal{O}^-_1$.  
Similarly, choose $\alpha_2$ small enough so that $A_{\alpha_2}( f(\mathcal{O}^+)) \Subset \mathcal{O}^+$ and $A_{\alpha_2}( f(\mathcal{O}^+_1)) \Subset \mathcal{O}^+_1$. Now consider the map $g$ defined as:
$$ g := A_{\alpha_2} \circ f \circ A_{\alpha_1}.$$
 The following properties are then satisfied:
\begin{itemize}
 \item $g^{-1}(\mathcal{O}^-)  \Subset \mathcal{O}^-$ and $g^{-1}(\mathcal{O}^-_1)  \Subset \mathcal{O}^-_1$ ;
\item  $g (\mathcal{O}^+)\Subset \mathcal{O}^+$ and $g (\mathcal{O}^+_1)\Subset \mathcal{O}^+_1$.
\end{itemize}
The example we have constructed is in the orbit of $f$ under the group $\mathrm{Aut}(\P^k)^2$ but it is of no concern since $f$ and $A_{\alpha_2}\circ f\circ A_{\alpha_2}^{-1}$ are conjugated.
Observe that the previous hypotheses are stable under small perturbations (that is conjugating with $c$ in a sufficiently small neighborhood of the identity in $\mathrm{Aut}(\P^k)$).  We deduce that:
\begin{lemme}\label{U,U_1...}
There exist analytic sets $E^-_{k-s}$, $E^-_{k-s-1}$ (resp. $E^+_{s}$ and  $E^+_{s-1}$) of dimension $k-s$, $k-s-1$, (resp. $s$ and $s-1$) and $\delta$-neighborhoods $\mathcal{O}^-$ and $\mathcal{O}^-_1$ (resp. $\mathcal{O}^+$ and $\mathcal{O}^+_1$)  of $E^-_{k-s}$ and $E^-_{k-s-1}$
(resp. $E^+_{s}$ and $E^+_{s-1}$) and an open set $\mathcal{Y}$ in the set of parameters such that for $c\in \mathcal{Y}$ and $f_c:= f\circ c$, we have:
\begin{itemize}
\item The sets $f_c^{-1}(E^-_{k-s})$ and $f_c^{-1}(E^-_{k-s-1})$ (resp. $f_c(E^+_{s})$ and $f_c(E^+_{s-1})$ ) are analytic sets of dimension $k-s$ and $k-s-1$ (resp. $s$ and $s-1$). 
\item $\overline{\mathcal{O}^-} \cap \overline{\mathcal{O}^+_1} =\varnothing$ and $\overline{\mathcal{O}^-_1} \cap \overline{\mathcal{O}^+} =\varnothing$.
\item  $f_c^{-1}(\mathcal{O}^-)  \Subset \mathcal{O}^-$ and $f_c^{-1}(\mathcal{O}^-_1)  \Subset \mathcal{O}^-_1$ (resp. $f_c(\mathcal{O}^+)\Subset \mathcal{O}^+$ and $f_c(\mathcal{O}^+_1)\Subset \mathcal{O}^+_1$).
\end{itemize}
\end{lemme}

We can now apply Proposition \ref{iterated_pull_back} to such an $f_c$ with $M=E^-_{k-s}$, $F=E^+_{s-1}$ (and $N=E^+_{s}$, $E=E^-_{k-s-1}$) and any $n$. It proves Theorem \ref{algebraic_stability_generic}.\\

We denote in what follows $d_q=\lambda_q(f)$ the generic dynamical degree in the orbit of $f$.

\section{Green currents in the generic case}
From now on,  we assume that the dynamical degree $d_s$ is strictly larger than $d_{s-1}$ (hence we have $1 < d_1 < \dots <d_s$).
Let $W := \mathrm{Aut}(\P^k)$. It is a Zariski dense open
set in the projective space $\widetilde{W}=\mathbb{P}^l$ where $l=(k+1)^2-1$. Let $c$ denote the homogeneous coordinate on
$\widetilde{W}$. When $c \in W$, we write $f_c$ instead of
$f\circ c$. We can extend the notation for $c \in \widetilde{W}$. Of
course, in that case $f_c$ is not a dominant meromorphic map and it might not be defined. For convenience, we denote $X:= \widetilde{W} \times \P^k$, it has complex dimension $(k+1)^2+k-1$.\\

Consider the rational map:
\begin{align*}
\widetilde{F}: X &\to X\\
(c,z) &\mapsto  (c,f_c(z)).
\end{align*}
Observe that $\widetilde{F}$ acts as the identity on $\widetilde{W}$. Let $\Pi_i$ denote the canonical projections of $X=\widetilde{W}\times\P^k$ to its factor for $i=1,2$. In $X$, let $\omega_i:=\Pi_i^*(\omega_{FS})$ be the pull-back
of the Fubini-Study form by the projection for $i=1,2$. That way, $\omega_1+\omega_2$ is a K\"ahler form on $X$. Let $\widetilde{\Gamma}$ be the graph of $\widetilde{F}$ in $X^2$ and let $P_i$, $i=1,2$, denote the projection from $X^2$ to its factors. We denote $\omega_{i,j}:= P_i^*(\omega_j)$.\\

Define $\widetilde{T}:=d^{-1}_s [\widetilde{\Gamma}]$ which is a positive closed current on $X\times X$ of bidegree $((k+1)^2+k-1,(k+1)^2+k-1)$.
 Let $\widetilde{\Omega}:=\sum_{a+b+c+d= (k+1)^2+k-1 } m_{a,b,c,d} \omega_{1,1}^a \wedge \omega_{1,2}^b \wedge \omega_{2,1}^c \wedge \omega_{2,2}^d$ 
be a smooth form cohomologous to $\widetilde{T}$. Such $\widetilde{\Omega}$ exists and is positive since $X$ is a product of projective spaces. One can then consider a negative quasi-potential $\widetilde{V}$ of $\widetilde{T}$ 
(that is $dd^c \widetilde{V}= \widetilde{T}-\widetilde{\Omega}$) given by Theorem 2.3.1 in \cite{DS6}.\\

We consider a sequence  $(\widetilde{T}_{m})$ of smooth positive closed currents such that:
 \begin{itemize}
\item each  $\widetilde{T}_{m}$ is cohomologous to $\widetilde{T}$;
\item $\widetilde{T}_{m} \to T$ in the sense of currents;
\item one can choose negative quasi-potentials $\widetilde{V}_{m}$ of $\widetilde{T}_m$ such that for all smooth positive closed currents $S$ of bidegree $((k+1)^2+k,(k+1)^2+k)$ in $X^2$, one has that $\langle  \widetilde{V}_m, S \rangle$ decreases to $\langle  \widetilde{V}, S \rangle$.
 \end{itemize}
Some remarks are in order here. Such sequence of currents has been explicitly constructed in \cite{DS6} where the authors restricts themselves to the case of $\P^k$ for simplicity. 
In order to construct $(\widetilde{T}_{m})$, one uses a convolution by a radial approximation of the identity in $\mathrm{Aut}(X^2)$ (more precisely, a poly-radial approximation). The last property is proved as in Proposition 3.1.6 in \cite{DS6}. Extending the formalism of super-potentials to $X^2$, we can extend that property to the case where $S$ is not smooth. Finally, that property implies the Hartogs' convergence of the sequence $\widetilde{T}_{m}$ to $\widetilde{T}$. We need some notations. For a current $R$ in $X$, we denote $\widetilde{L}(R):=(P_1)_*(P_2^*(R)\wedge \widetilde{T})$ and 
  $ \widetilde{\Lambda}(R):=(P_2)_*(P_1^*(R)\wedge \widetilde{T})$ (in the cases where these currents make sense) and
$ \widetilde{L}_{m}(R):=(P_1)_*(P_2^*(R)\wedge \widetilde{T}_{m})$ and 
  $ \widetilde{\Lambda_{m}}(R):=(P_2)_*(P_1^*(R)\wedge \widetilde{T}_{m})$. We have the lemma:
\begin{lemme}\label{decreasing_superpotentials}
Let $\Omega^q$ be smooth positive closed current of bidegree $(q,q)$ in $X$. Then one can choose negative quasi-potentials $\widetilde{U}$ and $\widetilde{U}_{m}$ of $\widetilde{L}(\Omega^q)$ and  $\widetilde{L}_{m}(\Omega^q)$ and $\widetilde{U}'$ and $\widetilde{U}'_{m}$ of  $\widetilde{\Lambda}(\Omega^q)$ and  $\widetilde{\Lambda}_{m}(\Omega^{q})$ such that 
\begin{itemize}
\item for all positive smooth forms $S$ of bidegree $(l-q+1,l-q+1)$ on $X$, we have that $\langle  \widetilde{U}_{m}, S \rangle$ decreases to $\langle  \widetilde{U}, S \rangle$. 
\item for all positive smooth forms $R$ of bidegree $(l-q+1,l-q+1)$ on $X$, we have that $\langle  \widetilde{U}'_{m}, R \rangle$ decreases to $\langle  \widetilde{U}', R \rangle$. 
\end{itemize}
\end{lemme}
\noindent \emph{Proof.} 
Let $\widetilde{U}$ and $\widetilde{U}_{m}$ be the negative quasi-potentials defined by 
$$\widetilde{U}:= (P_1)_*( P_2^*(\Omega^q)\wedge \widetilde{V}) \ \mathrm{and} \ \widetilde{U}_{m}:= (P_1)_*( P_2^*(\Omega^q)\wedge  \widetilde{V}_{m}).$$ 
Since $dd^c$ commutes with pull-back and push-forward we have indeed that  $dd^c \widetilde{U}= \widetilde{L}(\Omega^q)- (P_1)_*( P_2^*(\Omega^q)\wedge \widetilde{\Omega})$ 
and $dd^c \widetilde{U}_{m}= \widetilde{L}_{m}(\Omega^q)- (P_1)_*( P_2^*(\Omega^q)\wedge \widetilde{\Omega})$. 
The first part of the lemma follows from the choice of  $\widetilde{T}_{m}$ and $\widetilde{V}_m$. The proof of the second point is the same.  \hfill $\Box$ \hfill   \\

We shall now change the choice of the quasi-potentials. One of the interests of the theory super-potentials lies in the fact that although it is defined at some point using quasi-potentials, it does not depend on the choice of the quasi-potentials (up to a normalization). So we choose instead   $\widetilde{U}$ and $\widetilde{U}_{m}$ equal to the Green quasi-potentials of $\widetilde{L}(\Omega^q)$ and  $\widetilde{L}_{m}(\Omega^q)$ (see \cite[Theorem 2.3.1]{DS6}).

Recall that $\C_s$ is the set of normalized positive closed currents of bidegree $(s,s)$ in $\P^k$. Consider a smooth $\Omega_s \in \C_s$ (we will take more specific $\Omega_s$ later on). Consider $\Omega_s:=\Pi_2^*(\Omega_s)$. Then  $\Omega_s$ is a smooth positive closed current of bidegree $(s,s)$ in $X$. We apply the above lemma to $\Omega_s$. That gives negative quasi-potentials $\widetilde{U}$ and $\widetilde{U}_{m}$ of $\widetilde{L}(\Omega_s)$ and  $\widetilde{L}_m(\Omega_s)$. We let $\mathcal{U}$ and $\mathcal{U}_{m}$ be the associated super-potentials.
For $c \in W$, recall that $L_c$ and $\Lambda_c$ are the  corresponding normalized pull-back and push-forward operators.\\

We will need some tools on slicing theory and on convergence of DSH functions. Recall some facts on slicing first (see \cite{Fed} or \cite{DS12}[p. 483]). Let $\lambda_W$ be the standard volume form on $W$.
Let $\psi(c')$ be a positive smooth function with compact support in a chart of $W$ containing $c$
such that $\int\psi\lambda_W=1$. Define $\psi_\epsilon(c'):=\epsilon^{-2l}\psi(\epsilon^{-1} c')$ and 
$\psi_{c,\epsilon}(c'):=\psi_\epsilon(c'-c)$. The measures $\psi_{c,\epsilon}\lambda_W$
approximate the Dirac mass at $c$. Let $T$ be a current on $X$. For every smooth test form $\Psi$
on $X$
one defines the slice of $T$ at $c\in W$ as
\begin{equation} \label{eq_slice}
\langle T,\Pi_1, c \rangle (\Psi):=\lim_{\epsilon\rightarrow 0}
\langle  T \wedge \Pi_1^*(\psi_{c,\epsilon}\lambda_W),\Psi\rangle
\end{equation}
when $\langle T,\Pi_1 ,c\rangle$ exists.
This property holds for all choice of $\psi$.
Conversely, when the previous limit exists and is independent of $\psi$, 
it defines the current $\langle T,\Pi_1 ,c\rangle$  and one says 
that $\langle T,\Pi_1 ,c\rangle$ {\it is well defined} (similarly for  $\langle \widetilde{U}, \Pi_1, c \rangle$). \\

Restating results of Dinh and Sibony, we have the lemma :
\begin{lemme}\label{slice_n=1} Let $T$ be a positive closed current of bidegree $(q,q)$ on $X$ of super-potential $\U_{T}$, then outside a pluripolar set of $W$, one has that: 
\begin{itemize}
\item the slice $\langle T, \Pi_1, c \rangle$ is a well defined positive closed current on $\P^k$. 
\item the function $\gamma_{k-q+1}\mapsto \U_{T}(\Pi_2^*(\gamma_{k-q+1})\wedge[c'=c])$ is finite and  equal to a super-potential of $\langle T, \Pi_1, c \rangle$ on smooth forms.
\item if  $(T_m)$ is a sequence of positive closed currents on $X$ which converges in the Hartogs' sense to $T$, then the slices $\langle T_n, \Pi_1, c \rangle $ converges to $\langle T, \Pi_1, c \rangle $ in the Hartogs' sense.   
\item Finally, $\langle \widetilde{L}(\Omega_q), \Pi_1, c \rangle = \frac{d_q}{d_s}L_c(\Omega_q)$.  
\end{itemize} 
\end{lemme}
\noindent \emph{Proof.} The first point is proved in \cite{DS11}. For the second point, we consider a smooth form $\gamma_{k-q+1}$ of bidegree $(k-q+1,k-q+1)$. The quantity $\U_{T}(\Pi_2^*(\gamma_{k-q+1})\wedge[c'=c])$ is well defined 
by the theory of super-potentials (allowing the value $-\infty$). The set of $c$ where it is equal to $-\infty$ is pluripolar:
 else we can construct a probability measure $\mu$ on $W$ with bounded super-potential such that  $\U_{T}(\Pi_2^*(\gamma_{k-q+1})\wedge P_1^*(\mu))=-\infty$, a contradiction. Now, if $\U_{T}(\Pi_2^*(\gamma_{k-q+1})\wedge[c'=c])\neq -\infty$ for one smooth form $\gamma_{k-q+1}$, then it is true for any other smooth form as any smooth form is more H-regular than $\gamma_{k-q+1}$. 
 Observe now that for $c$ such that $\U_{T}(\Pi_2^*(\gamma_{k-q+1})\wedge[c'=c])$ is finite, we have that the quantity is equal to $\langle \langle U_{T}, \Pi_1, c \rangle,    \Pi_2^*(\gamma_{k-q+1}) \rangle$ (here, $U_{T}$ is a quasi-potential of $T$) : that follows from the definition of slicing. Hence,  $\mathcal{U}_{T} (P_2^*(\gamma_{k-q+1})\wedge[c'=c])$ defines a super-potential of $\langle T, \Pi_1, c \rangle $ on smooth forms. \\

The previous point means that $T$ and $[c'=c]$ are wedgeable. The third point follows  from Proposition 4.2.6 in \cite{DS6} as H-convergence is preserved by wedge-product.\\
 
For the last point,  the result is clear outside any neighborhood of $I(f_c)$ as $\widetilde{L}(\Omega_q)$ is continuous there. The result follows as the mass in a small neighborhood of $I(f_c)$ can be taken arbitrarily small since the mass of $L_c(\Omega_q)$ is $1$ (observe that a neighborhood of $I(f_c)$ is also a neighborhood of $I(f_{c'})$ for $c'$ close to $c$).  \hfill $\Box$ \hfill   \\

Recall some fact on dsh functions. We say that a function is quasi plurisubharmonic (qpsh for short) on $\widetilde{W}$ if is locally the difference of a plurisubharmonic function and a smooth function. 
  
We say that a measure is $PLB$ if the qpsh functions are integrable for that measure. Let $\mu$ be such a measure (any measure given by a smooth distribution for example).  We have the following lemma (see Proposition 2.4 in \cite{DS2}):
\begin{lemme}
 The family of qpsh functions in $\widetilde{W}$ such that
 $dd^c \psi \geq - \omega_1$ 
and one of the two following conditions:
 $$\max_{\widetilde{W}} \psi =0 \ \mathrm{or} \ \int \psi d\mu=0 $$
is bounded  in $L^1(\nu)$ and is bounded from above.
\end{lemme} 
We say that a function $u$ is dsh if can be written outside a pluripolar set as the difference of two qpsh functions. Then 
$dd^cu = T^+-T^-$ where $T^{\pm}$ are positive closed $(1,1)$-currents of same mass. 

For such a $u$, define
$$\|u\|_{DSH} := \left|\int u d\mu \right| + \min \|T^\pm\|$$
where the minimum is taken on all $T^\pm$ positive closed such that $dd^cu = T^+-T^-$. From a sequence of dsh functions uniformly bounded in $DSH$-norm,
one can extract a weakly converging subsequence (in the sense of currents with the norm of the limit bounded by the bound).

\begin{lemme}\label{merci_sib} Let $(g_n)$ be a bounded sequence of dsh functions on $W$. Then, we can extract a converging subsequence in $DSH$ that converges outside a pluripolar set.
 Then $g(x)$ is dsh and  $\|g\|_{DSH}\leq C$.
\end{lemme}
\noindent \emph{Proof.} Write $g_n= g^+_n-g^-_n$ where $g^\pm_n$ are qpsh functions such that
 $\|g^\pm_n\|_{L^1}$ is uniformly bounded (for the Lebesgues measure). Up to extracting,
 we can assume that $(g^{\pm}_n)$ converges outside a pluripolar set to $g^\pm$ (Proposition 3.9.4 in \cite{DS3}). \hfill $\Box$ \hfill   \\

\noindent{$\bullet$ \bf Action of $\widetilde{L}_{m}$ and  $\widetilde{L}$ on the cohomology.}
As $\widetilde{T}_{m}$ and $\widetilde{T}$ are cohomologous, $\widetilde{L}$ and  $\widetilde{L}_m$ coincide on the cohomology. We study the action of $\widetilde{L}$ on  $\omega_1^i\wedge\omega_2^{s-i}$. We have that $\widetilde{L}(\omega_1^i\wedge\omega_2^{s-i})$ is a positive closed form. Since $\widetilde{L}$ acts trivially on $\widetilde{W}$ we have that  $\widetilde{L}(\omega_1^i\wedge\omega_2^{s-i})= \omega_1^i\wedge\widetilde{L}(\omega_2^{s-i})$. 
 We can write it in cohomology (that is up to a $dd^c$-exact form):
 $$\widetilde{L}(\omega_1^i\wedge\omega_2^{s-i})=\omega_1^i\wedge \sum_{j=0}^{j=s-i}C_{j,s-i}\omega_1^j\wedge \omega_2^{s-i-j},$$
where the $C_i$ are non negative numbers (since $X$ is a product of projective space). We claim that
$$ C_{0,s-i}=\frac{d_{s-i}}{d_s}.$$
Indeed, Lemma \ref{slice_n=1} implies that for $c$ generic we have $\langle \widetilde{L}(\omega_2^{s-i}), \Pi_2, c\rangle =\frac{1}{d_s} f_c^*(\omega_2^{s-i})$ and then $C_{0,s-i}$ is just $\widetilde{L}(\omega_2^{s-i})$ evaluated at $[c'=c]\wedge \omega_2^{k-s}$.  
The matrix $\widetilde{M}$ of $\widetilde{L}$ on the basis $(\omega_1^i\wedge \omega_2^{s-i})$ is then the matrix with non negative coefficients:
$$ 
\left( \begin{array}{lllll}
1          & 0                         &  \cdots                &  \cdots          & 0 \\
\star      & \frac{d_{s-1}}{d_s}      & \ddots                 &  \               & \vdots  \\
\vdots     & \ddots                    & \frac{d_{s-2}}{d_s}   &  \ddots          & \vdots \\
\vdots     &                           & \ddots                 & \ddots                & 0 \\ 
\star      &                    \cdots &     \cdots             & \star            &  \frac{d_0}{d_s}
\end{array} \right).
$$
Using Perron-Frob\'enius theorem gives an eigenvector associated to the eigenvalue $1$ with positive coefficients. In other words, one can choose a positive closed form $\Omega^s$ such that $\widetilde{L}(\Omega^s)= \Omega^s+dd^c U_s$ and $\widetilde{L}_m(\Omega^s)= \Omega^s+dd^c U_{s,m}$. Renormalizing, the form $\Omega^s$ can be written as $\omega_2^s+ \sum_{i\geq 1} a_i \omega_1^i\wedge \omega_2^{s-i}$. Taking $\Omega_s$ (which is cohomologous to $\omega_2^s$) and $\Omega'_s$ any smooth form cohomologous to $\sum_{i\geq 1} a_i \omega_1^i\wedge \omega_2^{s-i}$  with no component of bidegree higher than $(s-1,s-1)$ in the $z$ variable ($z$ is the dynamical variable, i.e. the coordinate on the $\P^k$ factor), one can choose instead:
$$ \Omega^s:=\Omega_s+ \Omega'_s.$$

Similarly, the action of $\widetilde{\Lambda}$ on the basis  $(\omega_1^i\wedge \omega_2^{k-s-i})_{i=0..k-s}$  is given by the matrix:
$$ 
\left( \begin{array}{lllll}
1          & 0                         &  \cdots                &  \cdots          & 0 \\
\star      & \frac{d_{s+1}}{d_s}       & \ddots                 &  \               & \vdots  \\
\vdots     & \ddots                    & \frac{d_{s+2}}{d_s}   &  \ddots          & \vdots \\
\vdots     &                           & \ddots                 & \ddots                & 0 \\ 
\star      &                    \cdots &     \cdots             & \star            &  \frac{d_k}{d_s}
\end{array} \right).
$$
In order to see it, one can work with the dual basis $(\omega_1^{l-i}\wedge \omega_2^{s+i})_{i=0..k-s}$ and use that the adjoint of  $\widetilde{\Lambda}$ is  $\widetilde{L}$. \\

Since $\widetilde{L}(\Omega^s)=\widetilde{L}(\Omega_s)+ \widetilde{L}(\Omega'_s)$, we have $U_{s}=\widetilde{U}+\widetilde{U}'_s$ where $\widetilde{U}'_s$ is the Green quasi-potential of $\widetilde{L}(\Omega'_s)$ (recall that $\widetilde{U}$ is a quasi-potential of $\widetilde{L}(\Omega_s)$). Then $\widetilde{U}'_s$ is a form with no components of bidegree higher than $(s-2,s-2)$ in the $z$ variable. So bidegree arguments imply that $\langle U_{s}, \Pi_1, c \rangle=\langle \widetilde{U}, \Pi_1, c \rangle$ defines a super-potential $\U_{ L_c(\Omega_s)}$ of $L_c(\Omega_s)$. \\

\noindent{$\bullet$ \bf Construction of a function that tests the convergence of the Green current.}

We fix $N\in \mathbb{N}$. Let $m=(m_1,m_2,\dots, m_N)\in \mathbb{N}^N$.
We define:
$$g^m_N:= (\Pi_1)_*( \sum_{j=1}^{N} \widetilde{L}_{m_N} \dots\widetilde{L}_{m_{j+1}}(U_{s,m_j}) \wedge  \widetilde{\beta_{k-s+1}})
$$
where $\widetilde{\beta_{k-s+1}}= \Pi_2^{*}(\beta_{k-s+1})$ and $\beta_{k-s+1}$ is a positive closed current of bidegree $(k-s+1,k-s+1)$ 
in $\P^k$ and $U_{s,m_j}$ is a quasi-potential of $\widetilde{L}_{m_j}(\Omega^s)$.

Our aim is to prove the following proposition:
\begin{proposition}\label{dsh_approx}
\begin{enumerate}
\item There exist positive closed currents $T^\pm_{n,m}$ on $\widetilde{W}$ and a constant $C$ independent of $n$ and $m$ such that:
$$dd^c g^m_n =T^+_{n,m}-T^-_{n,m} $$
with $\| T^\pm_{n,m}\| \leq C$. 
\item Letting $m_N \to \infty, \ \dots, \ m_1\to \infty$ in that order, we have that the functions  $(g^m_n)$ converge outside a pluripolar set to the function $g_n$ defined for $c \in \cap_{n \in \Nn}\mathcal{Z}_{n,s}$ by:
 \begin{align}\label{conv_approx_test}
g_{n}(c)= \sum_{j\leq n} \left( \frac{d_{s-1}}{d_s}\right)^j \U_{ L_c(\Omega_s)}( \Lambda_c^{j}(\beta_{k-s+1})),
\end{align}
where $\U_{ L_c(\Omega_s)}$ is the super-potential of $L_c(\Omega_s)$ given on smooth forms by the quasi-potential $\langle U_s,\Pi_1, c\rangle$. 
\end{enumerate}
\end{proposition}

\noindent $\bullet$ {\bf Computation of $dd^c g^m_N$}
As every object in the definition of $g^m_N$ is smooth, its $dd^c$ is well defined (that is the very reason we introduced the regularization of the graph). Furthermore, $dd^c$ commutes with pull-back
and push-forward
\begin{align*}
dd^c g_N^m &=(\Pi_1)_*( \sum_{j=1}^{N} \widetilde{L}_{m_N} \dots\widetilde{L}_{m_{j+1}}(dd^c U_{s,m_j}) \wedge  \widetilde{\beta_{k-s+1}})\\
           &=(\Pi_1)_*( \sum_{j=1}^{N} \widetilde{L}_{m_N} \dots\widetilde{L}_{m_{j+1}}(\widetilde{L}_{m_j}(\Omega^{s})- \Omega^s) \wedge  \widetilde{\beta_{k-s+1}})\\
            &= (\Pi_1)_*(  \widetilde{L}_{m_N} \dots\widetilde{L}_{m_{1}}(\Omega^{s})\wedge   \widetilde{\beta_{k-s+1}}) - (\Pi_1)_*(  \Omega^{s}\wedge   \widetilde{\beta_{k-s+1}}).
\end{align*}
Hence, writing:
 \begin{align*}
T^+_{N,m} &:= (\Pi_1)_*(  \widetilde{L}_{m_N} \dots\widetilde{L}_{m_{1}}(\Omega^{s})\wedge   \widetilde{\beta_{k-s+1}})\\
T^-_{N,m} &:=(\Pi_1)_*(  \Omega^{s}\wedge   \widetilde{\beta_{k-s+1}}),
\end{align*}
we can write $dd^cg_N^m$ as the difference of two positive closed currents $T^+_{N,m}-T^-_{N,m}$. Since $T^-_{N,m}$ does not depends on $N$ and $m$, its mass is constant. 
Now $T^+_{N,m}$ has the same mass as  $T^-_{N,m}$ since they are cohomologous. 
That proves the first point of Proposition \ref{dsh_approx}.\\

 \noindent $\bullet$ {\bf Proof of the convergence of $g^m_N(c)$}\\
We can write that:
\begin{align*}
g_N^m(c) &=  \sum_{j=1}^{N} \widetilde{L}_{m_N} \dots\widetilde{L}_{m_{j+1}}(U_{s,m_j}) \wedge  \widetilde{\beta_{k-s+1}}\wedge [c'=c] \\
        &= \sum_{j=1}^{N} \langle \widetilde{U}_{s,m_j} , \widetilde{\Lambda}_{m_{j+1}}\dots \widetilde{\Lambda}_{m_{N}}  ( \widetilde{\beta_{k-s+1}} \wedge [c'=c]) \rangle.
\end{align*}
 Letting $m_N \to \infty$, we have that, for $c \in \cap_{n \in \Nn}\mathcal{Z}_{n,s}$, $\widetilde{\Lambda}_{m_{N}}  ( \widetilde{\beta_{k-s+1}} \wedge [c'=c])$ converges in the sense of currents to $\Pi_2^*(\frac{d_{s-1}}{d_s}\Lambda_c(\beta_{k-s+1}))\wedge[c'=c]$ (we can prove the convergence in the Hartogs' sense but we do not need it). Hence 
\begin{align*}\langle \widetilde{U}_{s,m_j} , \widetilde{\Lambda}_{m_{j+1}}\dots \widetilde{\Lambda}_{m_{N}}  ( \widetilde{\beta_{k-s+1}} \wedge [c'=c]) \rangle \to \\
 \langle \widetilde{U}_{s,m_j} , \widetilde{\Lambda}_{m_{j+1}}\dots \widetilde{\Lambda}_{m_{N-1}}  (\Pi_2^*(\frac{d_{s-1}}{d_s}\Lambda_c(\beta_{k-s+1}))\wedge[c'=c]) \rangle.
\end{align*}
  We let $m_{N-1}, \ \dots, \  m_{j+1}$ go to $\infty$ and we have that the previous quantity converges to  $\langle \widetilde{U}_{s,m_j} , \left( \frac{d_{s-1}}{d_s}\right)^j  \Pi_2^*( \Lambda_c^{j}(\beta_{k-s+1}))\wedge[c'=c] \rangle$ (at each step, all the objects but one are smooth so the convergence is clear). Now we let $m_j\to \infty$, Hartogs' convergence of $\widetilde{L}_{m_j}(\Omega^s)$ implies that the previous quantity converges to  $\langle U_{s} , \left( \frac{d_{s-1}}{d_s}\right)^j  \Pi_2^*( \Lambda_c^{j}(\beta_{k-s+1}))\wedge[c'=c] \rangle$. Thanks to the remark at the end of the paragraph where we computed the action on the cohomology, this can be rewritten as  $\left( \frac{d_{s-1}}{d_s}\right)^j \U_{ L_c(\Omega_s)}( \Lambda_c^{j}(\beta_{k-s+1}))$ . That proves Proposition \ref{dsh_approx}. \\

\noindent $\bullet$ {\bf Construction of the Green current for an open set in the space of parameters.}

We show now that the $g_N^m$ are uniformly bounded  in $\mathcal{Y}$. See Lemma \ref{U,U_1...} for the notations. We shall take for that specific $\Omega_s$, $\Omega_s'$ and $\beta_{k-s+1}$. Let $\Omega_s$ be a smooth positive closed current in $\mathcal{C}_s$ with support in 
$\mathcal{O}^-$. Let $\Omega'_s$ be a smooth positive closed current with support disjoint from $\mathcal{Y}$ (that can easily be done by choosing instead of $\omega_1$ a smooth approximation of a hyperplane not meeting $\mathcal{Y}$, restricting $\mathcal{Y}$ if necessary). Let $\beta_{k-s+1}$ be a smooth positive closed current in $\mathcal{C}_{k-s-1}$ with support in $\mathcal{O}^+_1$. 
Observe that by construction of $\widetilde{L}_m$ and $\widetilde{\Lambda}_m$, we have that for $m$ large enough $\mathrm{Supp}(\widetilde{L}_m(\Omega_s))$,  $\mathrm{Supp}(\widetilde{L}_m(\Omega'_s))$ and $\mathrm{Supp}(\widetilde{\Lambda}_m(\widetilde{\beta}_{k-s+1}))$ are close to $\mathrm{Supp}(\widetilde{L}(\Omega_s))$,  $\mathrm{Supp}(\widetilde{L}(\Omega'_s))$ and $\mathrm{Supp}(\widetilde{\Lambda}(\widetilde{\beta}_{k-s+1}))$. 

In particular for $c\in\mathcal{Y}$, we have that $\widetilde{L}_{m_j}(\Omega^s)=\widetilde{L}_{m_j}(\Omega_s)$ and it has support in $\mathcal{O}^-\times \mathcal{Y}$.
 Lemma 2.3.5 in \cite{DS6} implies that there is a constant $C>0$ (that does not depend on $m_j$) such that
$$\| U_{s,m_j}\|_{\C^1(\mathcal{O}^+_1\times \mathcal{Y})}\leq  C.$$
Slicing implies that: 
$$\|\langle U_{s, m_j}, \Pi_1, c\rangle \|_{\C^1(\mathcal{O}^+_1)}\leq  C,$$
for $c\in \mathcal{Y}$.

Since $\widetilde{\Lambda}_{m_{j+1}}\dots \widetilde{\Lambda}_{m_{N}}  ( \widetilde{\beta_{k-s+1}} \wedge [c'=c])$ has support in $\mathcal{O}^+_1$ and mass $(\frac{d_{s-1}}{d_s})^n$, we deduce that there exists a constant $C_0$ independent of $m$ and $n$ (providing that $m$ is large enough with respect to $n$) such that $g_n^m$ is uniformly bounded by $C_0$ for $c$ in $\mathcal{Y}$. \\

\noindent $\bullet$ {\bf Construction of the Green current outside a pluripolar set.}

Take $\mu$ a smooth measure with support in $\mathcal{Y}$. Such $\mu$ is PLB and is the one we use to define the DSH-norm. Then $n$ being fixed, we have that the sequence of functions $g_n^m$ is uniformly bounded in DSH, we can assume that it converges (in DSH). In particular, its limit $g'_n$ is DSH with $\|g'_n\|_{DSH}\leq C$ by Proposition \ref{dsh_approx} and $g'_n=g_n$ by Lemma \ref{merci_sib}. In particular, the sequence $(g_n)$ is uniformly bounded in $DSH$. Since the sequence of (non positive) functions $g_n$ is decreasing (and well defined outside a pluripolar set), we have that it converges for $c$ outside a pluripolar set to $g(c)$ in $\mathbb{R}^- \cup \{\infty\}$. On the other hand, we can extract a weakly converging sequence in $DSH$ to a limit $g'$. Extracting if necessary, we can assume that the convergence holds outside a pluripolar set by Lemma \ref{merci_sib}. In particular, $g=g'$ outside a pluripolar set. Hence, $g$ is finite outside a pluripolar set (removing if necessary the pluripolar set $(\cap_{n \in  \Nn} \mathcal{Z}_{n,s})^c$, we assume from now on that this pluripolar set contains it). \\ 

The sum
 \begin{align*}
g_{n}= \sum_{j\leq n} \left( \frac{d_{s-1}}{d_s}\right)^j \U_{ L_c(\Omega_s)}( \Lambda_c^{j}(.)),
\end{align*}
defines a super-potentials of $L_c^n(\Omega_s)$ by Theorem \ref{algebraic_stability_generic}. In here, the function $g_n$ is extended in addition to the parameter $c$ to a second argument, namely the input current $\beta_{k-s+1}$. One of the key points of super-potential theory, is that the finiteness of $g_n$ at $\beta_{k-s+1}$ implies the finiteness of $g_n$ at any current more H-regular than $\beta_{k-s+1}$ and in particular for all smooth forms. The sequence is decreasing and outside a pluripolar set, it does not converge to $-\infty$.
Outside that set, the convergence of the sequence implies the convergence in the Hartogs' sense of the sequence of currents $(L^n_c(\Omega_s))$
(see Corollary 3.2.7 in \cite{DS6}). We denote its limit by $T^+_{s,c}$ that we call the \emph{Green current of order $s$ of $f_c$}.
 Observe that the convergence of $(L^n_c(\Omega_s))$ in the Hartogs' sense to $T^+_{s,c}$ implies the convergence of $(L^n_c(\Theta_s))$ in the Hartogs' sense to $T^+_{s,c}$ for any other smooth form $\Theta_s \in \C_s$ (that is because any smooth form is more H-regular than any other current).  \\

The current $T^+_{s,c}$ is $f_c^*$-invariant (in the sense of super-potentials) since the convergence of the series giving $g_n(c)$ (see (\ref{conv_approx_test})) implies the convergence of:
 \begin{align*}
 \sum_{2\leq i} \left( \frac{d_{s-1}}{d_s}\right)^i \U_{ L_c(\Omega_s)}\circ \Lambda_c^{i}(\beta_{k-s+1}).
\end{align*}
 Factorizing, we get that $ \sum_{1\leq i} \left( \frac{d_{s-1}}{d_s}\right)^i \U_{ L_c(\Omega_s)}\circ \Lambda_c^{i}(\Lambda_c(\beta_{k-s+1}))$ converges hence a super-potential 
of $T^+_{s,c}$ is finite at $\Lambda_c(\beta_{k-s+1})$. That means that $T^+_{s,c}$ is $f_c^*$-admissible (see Definition 5.1.4 in \cite{DS6}) and $L_c(T^+_{s,c})$ is well defined. Now, $L^{n+1}_c(\Omega_s)= L_c L_c^n( \Omega_s)$ converges in the Hartogs' sense to $T^+_{s,c}$ and Proposition 5.1.8 in \cite{DS6} implies that it also converges to $L_c(T^+_{s,c})$. Thus $T^+_{s,c}=L_c(T^+_{s,c})$ and $T^+_{s,c}$ is $f_c^*$-invariant. In particular, we have proved:
\begin{theorem}
There exists a pluripolar set $\mathcal{P}$ of $W$ such that for any $c\notin W$
 for any smooth form $\Omega_s \in \C_s$, the sequence of currents $L^n_c(\Omega_s)$
 converges in the Hartogs' sense to the Green current $T^+_{s,c}$ which is $f_c^*$-invariant.
\end{theorem}

Now assume that furthermore, $s$ is such that $d_s>d_{s+1}>\dots > d_k$ so $d_s$ is the 
highest degree. In other words, we have that generic maps in the orbit of $f$ are cohomologically hyperbolic. Doing the same thing for $\Lambda_c$, we obtain:
\begin{theorem}\label{convergence_pp_forward}
There exists a pluripolar set $\mathcal{P}$ of $W$ such that for any $c\notin W$
 for any smooth form $\Omega_{k-s} \in \C_{k-s}$, the sequence of currents $\Lambda^n_c(\Omega_{k-s})$
 converges in the Hartogs' sense to the Green current $T^-_{s,c}$ which is $(f_c)_*$-invariant.
\end{theorem}

 \noindent $\bullet$ {\bf Wedge product of $T^+_{s,c}$ and $T^-_{s,c}$ outside a pluripolar set}\\
We now prove:
\begin{proposition}\label{intersection} 
Outside a pluripolar set, the currents $T^+_{s,c}$ and $T^-_{s,c}$ are wedgeable so the probability measure 
$T^+_{s,c}\wedge T^-_{s,c}$ is well defined.
\end{proposition}
\noindent\emph{Proof.} Recall that $T^+_{s,c}$ and $T^-_{s,c}$ are wedgeable if a super-potential of $T^+_{s,c}$ is finite at $\Omega_1 \wedge T^-_{s,c}$ for one smooth form $\Omega_1 \in \C_s$. Let $\beta_{k-s} \in \C_{k-s}$ be a smooth form.  We will choose particular $\Omega_1$ and $\beta_{k-s}$ later. 
Consider the lemma:
\begin{lemme}\label{pourg'}
The sequence of functions
$$ g'_{n,m}(c)= \sum_{i\leq n-1} \left( \frac{d_{s-1}}{d_s}\right)^i \U_{ L_c(\Omega_s)}(\Lambda_c^{i}(\Omega_1 \wedge \Lambda_c^{m}\beta_{k-s})) $$
is a sequence of DSH functions uniformly bounded in $n$ and $m$ for the DSH norm that converges  outside a pluripolar  set when $n\to\infty$ to $\U_{T^+_{s,c}}(\Omega_1 \wedge \Lambda_c^{m}\beta_{k-s})$.
\end{lemme}
Assume the lemma is proved, then we have that $\U_{T^+_{s,c}}(\Omega_1 \wedge \Lambda_c^{m}\beta_{k-s})$ defines a bounded sequence of $DSH$ functions. When $m\to \infty$, it converges outside a pluripolar set to $\U_{T^+_{s,c}}(\Omega_1 \wedge T^-_{s,c})$ since $\Lambda_c^{m}\beta_{k-s}$ converges to $T^-_{s,c}$ in the Hartogs' sense. Hence $\U_{T^+_{s,c}}(\Omega_1 \wedge T^-_{s,c})\neq-\infty$ outside a pluripolar set and the proposition is proved.\hfill $\Box$ \hfill   \\

\noindent \emph{Proof of the lemma.} In order to control the DSH norm of $g'_{n,m}$, we need to compute its $dd^c$. That is done exactly as in the proof of Proposition \ref{dsh_approx} replacing $\widetilde{L}$ and $\widetilde{\Lambda}$ by their smooth approximations in order to deal with smooth objects and  using  $\Omega_1 \wedge \Lambda_c^{n} (\beta^{k-s})$ instead of $\beta_{k-s+1}$ (where $\beta^{k-s}=\beta_{k-s}+\dots$ is the eigenvector of $\widetilde{\Lambda}$ associated to 1). So, all there is left is to construct a PLB measure 
$\mu$ for which $\|g'_{n,m}\|_{L^1(\mu)}$ is uniformly bounded. As in the previous section, that will be achieved
by constructing an example stable by pertubations for which
we have uniform estimates in the convergence of $g'_{n,m}$. \\ 

We use the notations and results of Lemma \ref{U,U_1...}. We consider parameters $c\in \mathcal{Y}$.
 As in the previous paragraph, we take $\Omega_s\in \C_s$ a smooth current with support in $\mathcal{O}^-$. We take $\beta_{k-s}\in \C_{k-s}$ any smooth form with support in $\mathcal{O}^+$. In particular, $\Lambda^m_c(\beta_{k-s})$ has support in  $\mathcal{O}^+$ for all $m$.

Let $H$ be the hyperplane spanned by $E^-_{s-1}$ and $E^-_{k-s-1}$. Let $\Omega_1$ be a smooth element of $\C_1$ with support in a small neighborhood of $H$. Choosing that neighborhood small enough, we have that  $\Omega_1 \wedge \Lambda^m_c(\beta_{k-s})$ is a probability measure with support in $\mathcal{O}^+_1$.  In particular, for $c\in \mathcal{Y}$, $\U_{ L_c(\Omega_s)}(\Lambda_c^{i}(\Omega_1 \wedge \Lambda^m_c(\beta_{k-s})))=\langle U_{ L_c(\Omega_s)}  ,\Lambda_c^{i}(\Omega_1 \wedge \Lambda^m_c(\beta_{k-s}))   \rangle$ as $U_{ L_c(\Omega_s)}$ is smooth on the support of $\Lambda_c^{i}(\Omega_1 \wedge \Lambda^m_c(\beta_{k-s}))$. For $c\in \mathcal{Y}$, we have that:
$$\|\langle U_{s}, \Pi_1, c\rangle \|_{\C^1(\mathcal{O}^+_1)}\leq  C.$$
Hence:
$$\left| \left( \frac{d_{s-1}}{d_s}\right)^i \U_{ L_c(\Omega_s)}(\Lambda_c^{i}(\Omega_1 \wedge \Lambda^m_c(\beta_{k-s})))\right| \leq C\left(\frac{d_{s-1}}{d_s}\right)^i.$$ 
That implies that $|\U_{L_c^n(\Omega_s)}(\beta_{k-s+1})|$ is uniformly bounded
by a constant $C_0$ in $\mathcal{Y}$ where $C_0$ does not depend on $n, m$. Again, we take for $\mu$ any smooth measure with support in $\mathcal{Y}$.  \hfill $\Box$ \hfill   \\

In particular, we have proved points 2 and 3 in Theorem \ref{principal}.

\begin{Remark} \normalfont 
\begin{enumerate}
\item Hartogs' regularity implies that, for $c$ generic,  $L_c^n(\Omega_s)$ and $\Lambda_c^m(\Omega_{k-s})$ are wedgeable for any smooth $\Omega_s$ and $\Omega_{k-s}$ and $n, m$.
\item It does not follows from the proposition that the measure $T^+_{s,c}\wedge T^-_{s,c}$ is invariant. 
Indeed, we have not proved that it does not charge $I(f_c)$.
If not, such measure would hold little interest for the dynamics of $f_c$. So, in the next section, we will show a (stronger) property of invariance (namely the quasi-potential of $T^+_{1,c}$, the Green current of order $1$ is integrable with respect to $T^+_{s,c}\wedge T^-_{s,c}$). 
\end{enumerate}
\end{Remark}

\section{Green measure in the generic case}
In that section, we assume again that $d_s$ is the largest (generic) dynamical degree. Our purpose is to prove the following which will give point 4 in Theorem \ref{principal} :
\begin{theorem}\label{1andeboulot}
Let $f$ be such that $\mathrm{dim}(I(f))=k-s-1$  or $I\subset H$ for a hyperplane $H$.
Then there exists a pluripolar set $\mathcal{P}$ of $W$ such that for any $c\notin W$ the Green currents $T^+_{s, c}$ and  $T^-_{s,c}$
are well defined, wedgeable. Furthermore, the measure $\nu_c:=  T^+_{s,c}\wedge T^-_{s,c}$ is an invariant probability measure that integrates $\log \mathrm{dist}( , I(f_c))$ of maximal entropy $\log d_s$.

The Lyapunov exponents $\chi_1 \geq \chi_2 \geq \dots \geq \chi_k$ of $\nu_c$ are well defined and we have the estimates:
\begin{align*}
\chi_1\geq\dots\geq\chi_s\geq \frac{1}{2}\log \frac{d_s}{d_{s-1}}>0\\
0>\frac{1}{2}\log \frac{d_{s+1}}{d_{s}} \geq \chi_{s+1}\geq\dots\geq\chi_k \geq -\infty. 
\end{align*}
In particular, the measure $\nu_c$ is hyperbolic. 
\end{theorem}

\noindent $\bullet$ {\bf Strategy of the proof.}
We shall construct the measure of maximal entropy using a theorem of De Th\'elin and  the
author (\cite{DV1}) :
\begin{theorem}\label{ENTROPY}
Consider the sequence of measures:
$$\nu_{c,n}:= \frac{1}{n} \sum_{i=0}^{n-1} (f_c^i)_{*}\left( \frac{(f_c^n)^* \omega^{s}
  \wedge \omega^{k-s}}{\lambda_l(f_c^n)} \right).$$
Assume that there exists a converging subsequence $\nu_{c,\psi(n)} \to \nu_c$ with:
$$(H) \mathrm{  :   }  \lim_{n \rightarrow + \infty}
\int \log d(x,I)  d \nu_{c,\psi(n)} (x) = \int \log d(x,I(f_c))  d \nu_c(x) > - \infty.$$
Then $\nu_c$ is an invariant measure 
of metric entropy $=\log d_s$.
\end{theorem}
Observe that in (\cite{DV1}), we define $\nu_c$ for $s$ not necessarily associated to the highest dynamical degree and then we only have that
$\nu_c$ is an invariant measure 
of metric entropy  $\geq \log d_s$. But in our case, the other inequality always stands by \cite{DS9}. \\

The estimates on the Lyapunov exponents follows from Corollary 3 in \cite{DT1}.

Observe that in that theorem, one requires that $\log \mathrm{dist}(x, \mathcal{A}) \in L^1(\nu_c)$ where  $ \mathcal{A}= C_{f_c} \cup I_{f_c}$ (recall that $C_{f_c}$ is the critical set of $f_c$). But in our case, we only have that $\log \mathrm{dist}(x, I_{f_c}) \in L^1(\nu_c)$. Despite that fact, one still has the hyperbolicity of the measure allowing the value $-\infty$ for the negative Lyapunov exponents. Indeed, the stable manifolds were obtained in \cite{DT1} by composing  forward graph transforms for $f^{-1}$ along $\nu_c$-generic orbits. In the non-integrable case, one can produce them by performing backward graph transforms for $f$ itself. Then, once the stable manifolds are constructed, volume estimates are obtained by the slicing
arguments of \cite{DT1} (we are very grateful to De Th\'elin for explaining that fact to us, one can also see \cite{Dup} where De Th\'elin's arguments are checked).

 Last, we do not claim that the Lyapunov exponents are constant (that is the case if $\nu_c$ is ergodic), but considering a ergodic decomposition of $\nu_c$, we have that almost all the measures appearing in the decomposition are ergodic (some could have mass on $I^+$, but only a set of $0$ measure since $\int \log d(x,I(f_c))  d \nu_c(x) > - \infty$).  Similarly, almost all the measures appearing in the decomposition are of maximal entropy $=\log d_s$ (because entropy is convex with respect to the measure and all of them are of entropy less than $\log d_s$). Finally, almost all the measures appearing in the decomposition integrates $\log d(x,I(f_c))$ (same reasons). Finally, we apply Corollary 3 in \cite{DT1} to each one of these generic measures of the decomposition.
 \\

In particular, Theorem \ref{1andeboulot} is proved if we can apply Theorem \ref{ENTROPY} for $c$ outside a pluripolar set.
We are going for that to follow the strategy of \cite[Proposition 3.4.16]{DV1}: 
one can apply Theorem \ref{ENTROPY} and obtain the following writing of $\nu_c$, providing we can prove the theorem: 
\begin{theorem}\label{constructionmeasure}
let $f$ such that $\mathrm{dim}(I(f))=k-s-1$  or $I\subset H$ for a hyperplane $H$.  
Outside a pluripolar set, the current  $T^+_{s,c}$ and $T^-_{s,c}$ are wedgeable. 
So the intersection $T^+_{s,c}\wedge T^-_{s,c}$ is a well defined probability measure $\nu_c$ and the quasi-potential of $T^+_{1,c}$, the Green current of order $1$, is integrable with respect to that measure.
\end{theorem} 
Assume the theorem is proved. Let us briefly explain how we can conclude. Since $L_c(\omega)$ is more H-regular than $T^+_{1,c}$, we also have that $\nu_c$ integrates a quasi-potential $U_{L_c(\Omega)}$ of $L_c(\omega)$.  Now, a quasi-potential of $L_c(\omega)$ has singularities in $\log d(x,I(f_c))$. Hartogs' regularities implies that $L^n_c(\omega^s)\wedge \Lambda_c^m(\omega^{k-s})$ is a well defined probability measure that integrates a quasi-potential  of $L_c(\omega)$. In particular, it does not charge $I(f_c)$ and it is $(f_c)_*$-admissible. Replacing $L_c^{n-1}(\omega^s)$ and $\Lambda_c^m(\omega^{k-s})$ by sequences of smooth currents converging in the Hartogs' sense, we prove that 
$$\Lambda_c(L^n_c(\omega^s)\wedge \Lambda_c^m(\omega^{k-s}))=L^{n-1}_c(\omega^s)\wedge \Lambda_c^{m+1}(\omega^{k-s}).$$ 
Again, continuity of the wedge product and $f_*$ for the H-convergence implies that $\nu_c$ is $f_*$-invariant and we can write $\nu_{c,n}$ as:
$$ \frac{1}{n} \sum_{i=0}^{n-1} L^{n-i}_c(\omega^s)\wedge \Lambda_c^i(\omega^{k-s}).$$
It follows that $\nu_{c,n}$ converges to $\nu_c$ in the Hartogs' sense and satisfies the condition (H). Then, we can apply Theorem \ref{ENTROPY}.\\

Observe also that the fact that $\nu_c$ integrates a quasi-potential of $L_c(\omega)$ is equivalent to the fact that it integrates a quasi-potential of $T^+_{1,c}$. Indeed, if $\nu_c$ integrates a quasi-potential of $L_c(\omega)$ it is $(f_c)_*$-invariant (see just above). A simple recurrence shows that it integrates $(f_c^*)^nU_{L_c(\Omega)}$ and 
$$\langle \nu_c , U_{L_c(\Omega)}\rangle =\langle \nu_c , (f_c^*)^nU_{L_c(\Omega)}\rangle. $$
The result follows by monotone convergence as a quasi-potential of $T^+_{1,c}$ is given by $\sum_{1\leq n} \frac{1}{d_1^n} (f_c^*)^nU_{L_c(\Omega)}$.\\

Proposition \ref{intersection} already states that $T^+_{s,c}$ and $T^-_{s,c}$ are wedgeable for $c$ generic. So we only need to prove
that the potential of the Green current of order $1$ is integrable with respect to $T^+_{s,c}\wedge T^-_{s,c}$ or, as it was observed in the above paragraph, that the potential of $L_c(\omega)$ is integrable with respect to $T^+_{s,c}\wedge T^-_{s,c}$.
We proceed as in the previous section. Let $\Omega_1 \in \C_1(\P^k)$, we consider $\Pi_2^*(\Omega_1)$ that we simply  denote by $\Omega_1$. 
Let $U_1$ be a quasi-potential of $\frac{1}{d_1}F^*(\Omega_1)$. Then outside a pluripolar set of $c$ one has that $ L_c(\Omega_1)=\langle \frac{1}{d_1}F^*(\Omega_1), \Pi_1, c\rangle $ and
 the slice $\langle U_1, \Pi_1, c\rangle $ is a quasi-potential 
of  $ L_c(\Omega_1)$ (in fact, that is true for all $c$). We denote by $\U_{ L_c(\Omega_1)}$ the associated super-potential.
Let $\beta_{k-s}\in \C_{k-s}$ (we will choose a more specific $\beta_{k-s}$ later on). Consider the lemma:
\begin{lemme}\label{pourk}
Let $f$ be such that $\mathrm{dim}(I(f))=k-s-1$  or $I\subset H$ for a hyperplane $H$.  
The sequence of functions
$$ k_{n}(c)=  \U_{ L_c(\Omega_1)}( L^n_c(\Omega_s)\wedge  \Lambda^n_c(\beta_{k-s}))$$ 
is a sequence of DSH functions uniformly bounded in $n$ for the DSH norm.
\end{lemme}
Assume the lemma is proved. From above, we have that outside a pluripolar set,  $L^n_c(\Omega_s)\wedge  \Lambda^n_c(\beta_{k-s})$
converges to $T^+_{s,c}\wedge T^-_{s,c}$ in the Hartogs' sense (Proposition 4.2.6 in \cite{DS6}). Hence,
$\U_{ L_c(\Omega_1)}( L^n_c(\Omega_s)\wedge  \Lambda^n_c(\beta_{k-s}))$ converges to 
$\U_{ L_c(\Omega_1)}( T^+_{s,c}\wedge T^-_{s,c})$ by continuity of the super-potential for 
the Hartogs' convergence (Remarks 3.2.4. in \cite{DS6}).
Then extracting weakly converging sequences in DSH to a limit $k$ and
using Lemma \ref{merci_sib}, we deduce that $k \neq -\infty$ outside a pluripolar set. As $k(c)=\langle U_{ L_c(\Omega_1)}, T^+_{s,c}\wedge T^-_{s,c}\rangle$ that implies Theorem \ref{constructionmeasure}. \\

In order to prove Lemma \ref{pourk}, we first have to control $dd^c k_n$. 
That is done exactly as above using the same techniques of approximation 
in the Hartogs' sense of the graph of the application $\widetilde{F}$. 
So, all there is left is to construct is the PLB measure 
$\mu$ on $\widetilde{W}$ such that 
$\|k_{n}(c)\|_{L^1(\mu)}$ are uniformly bounded. As in the previous section, that will be achieved
by constructing an example stable by pertubations for which
we have uniform estimates in the convergence of $k_n$.
We will first do that in the case where $\mathrm{dim}(I)=k-s-1$ and then when $I\subset H$ for some hyperplane $H$. \\ 

\noindent $\bullet$ {\bf Construction of an example stable by perturbations when $\mathrm{dim}(I)=k-s-1$.}

Recall that we constructed linear subspaces 
$E^+_s$, $E^+_{s-1}$, $E^-_{k-s-1}$, $E^-_{k-s}$ 
in Lemma \ref{U,U_1...}. We can assume that 
$I \cap \overline{\mathcal{O}^+} =\varnothing$ since $\mathcal{O}^+$ is a small neighborhood of a linear set of dimension $s$.
 As in Section \ref{Algebraic stability is dense}, we choose 
an element $A_{\alpha}$ in $\mathrm{Aut}(\P^k)$ such that 
\begin{itemize}
\item $A_{\alpha}^{-1}(f^{-1})(\mathcal{O}^-) \Subset \mathcal{O}^-$,
\item $A_{\alpha}^{-1}(f^{-1}(\mathcal{O}^-_1)) \Subset \mathcal{O}^-_1$,  
\item $A_{\alpha}( f(\mathcal{O}^+)) \Subset \mathcal{O}^+$, 
\item $A_{\alpha}( f(\mathcal{O}^+_1)) \Subset \mathcal{O}^+_1$. 
\end{itemize}
Consider the element $g$ defined as:
$$ g := A_{\alpha} \circ f \circ A_{\alpha}.$$
 Observe that $I(g)= A_{\alpha}^{-1}(I(f))$ hence we can assume (taking $\alpha$ large enough) that $I(g) \subset  \mathcal{O}_1^-$.  
 The following property are then satisfied:
\begin{itemize}
 \item $g^{-1}(\mathcal{O}^-)  \Subset \mathcal{O}^-$ and $g^{-1}(\mathcal{O}^-_1)  \Subset \mathcal{O}^-_1$ ;
\item  $g (\mathcal{O}^+)\Subset \mathcal{O}^+$ and $g (\mathcal{O}^+_1)\Subset \mathcal{O}^+_1$;
\item  $I(g) \subset  \mathcal{O}_1^-$.
\end{itemize}
Again, the example we have constructed is in the orbit of $f$ under the group $\mathrm{Aut}(\P^k)^2$ but that is of no concern since $f$ and $A_{\alpha}\circ f\circ A_{\alpha}^{-1}$ are conjugated.
Observe that the previous properties are stable under small perturbations, so we can find a smooth probability measure $\mu$ with support in 
$W$ such that the above conditions are satisfied for $g_c=g\circ c$ with $c\in \mathrm{Supp}(\mu)$. 

Now, in order to prove Lemma \ref{pourk}, we choose for $\Omega_1$ any smooth form in $\C_1$ (for example the Fubini-Study form). 
As before, we take for $\Omega_s$ a smooth form in 
$\C_s$ with support in $\mathcal{O}^-$ and for $\beta_{k-s}$
 a smooth form in $\C_{k-s}$ with support in 
$\mathcal{O}^+$. In particular, $L_n^c(\Omega_s)$ has support in $\mathcal{O}^-$ and $\Lambda_c^n(\beta_{k-s+1})$ has support in $\mathcal{O}^+$.
Thus $L^n_c(\Omega_s)\wedge  \Lambda^n_c(\beta_{k-s})$ is a probability measure with 
support in $\mathcal{O}^+\cap \mathcal{O}^-$. The super-potential $\U_{ L_c(\Omega_1)}$
is given by a quasi-potential $U_{ L_c(\Omega_1)}$. Lemma 2.3.5 in \cite{DS6} implies that there is a constant $C>0$ independent of $c$ such that
$$\|U_{ L_c(\Omega)}\|_{\C^1(\mathcal{O}^+_{s}\cap \mathcal{O}^-_{k-s})}\leq  C.$$
So arguing as above, we have that $k_n(c)$ is uniformly bounded for $c\in \mathrm{Supp}(\mu)$. That gives Lemma  \ref{pourk} in
the case where $\mathrm{dim}(I)=k-s-1$. \\

\noindent $\bullet$ {\bf Construction of an example stable by perturbations when $I$ is contained in a hyperplane.}

We modify the previous construction. Let $H$ denote a hyperplane such that $I \subset H$.
 Let $E^+_s$ and $E^-_{k-s}$ be (generic) linear subspaces of $\P^k$ of dimension $s$ and $k-s$. 
We consider $E_{s-1}^+:=  E_s^+ \cap H$ and $E_{k-s-1}^-:= E_{k-s}^-\cap H$. Then $E_{s-1}^+$ and $E_{k-s-1}^-$ are linear subspaces 
of dimension $s-1$ and $k-s-1$. 
We claim that we can assume:
\begin{itemize}
\item $E^+_s \cap E^-_{k-s} =\{p\}$ is reduced to a point and $H \cap \{p\} =\varnothing$ 
\item $E^+_s$  (resp. $E^-_{k-s}$) is $f_*$-compatible (resp. $f^*$-compatible)
\item $f(E^+_s) \cap E_{k-s-1}^-=\varnothing $ and $f^{-1}(E^-_{k-s}) \cap E_{s-1}^+=\varnothing $
\end{itemize}
We explain why the last point stands. It is generic (in the algebraic sense), hence we only need to show that it is not empty. For that we  can choose $E^-_{k-s}$ so that $f^{-1}(E^-_{k-s}) \cap H$ is of dimension $k-s-1$. Indeed the set of $Z\in E^-_{k-s}$ sent to $H$ by $f$ is a proper analytic set of $E^-_{k-s}$ and is then of dimension $\leq k-s-1$. We conclude using the first point of Lemma \ref{pull_back_generic_linear}. Since $\mathrm{dim}(f^{-1}(E^-_{k-s}) \cap H)+\mathrm{dim}(E^+_{s-1})= k-2 <\mathrm{dim}(H)$, we can assume that $f^{-1}(E^-_{k-s}) \cap E_{s-1}^+=\varnothing $. We  proceed similarly for $f(E^+_s) \cap E_{k-s-1}^-=\varnothing $.

We let $\mathcal{O}^+$ be a small neighborhood of $E^+_s$ and $\mathcal{O}^-$
be a small neighborhood of $E^-_{k-s}$. We can assume that $ H \cap \overline{\mathcal{O}^+\cap \mathcal{O}^- } =\varnothing$. 
We choose small neighborhoods $\mathcal{O}^+_1$ and $\mathcal{O}^-_1$ of $E_{s-1}^+$ and $E_{k-s-1}^-$. 
We can choose the homogeneous coordinates $[z_0: \dots: z_k]$ such that
\begin{align*}
E^-_{k-s}=\{z_{0}=\dots= z_{s-1}=0 \}&  \ \mathrm{and} \quad E^+_{s}=\{z_{s+1}=\dots= z_{k}=0 \} \\
E^-_{k-s-1}=\{z_{0}=\dots= z_{s-1}={z_s}=0 \}& \ \mathrm{and} \quad E^+_{s-1}=\{z_{s}=\dots= z_{k}=0 \}\\
H=\{z_s=0\}
\end{align*}
As in section \ref{Algebraic stability is dense} let $A_{\alpha}$ be the element of $\mathrm{Aut}(\P^k)$ given by 
\begin{align*} 
&A_{\alpha} ([z_0 : z_1 : \dots z_{s-1} : z_s: z_{s+1} :\dots: z_k ])=  \\
&\quad [ \alpha^{-1} z_0 :\alpha^{-1} z_1 : \dots \alpha^{-1} z_{s-1} : z_s: \alpha z_{s+1} :\dots: \alpha z_k ]. 
 \end{align*}  
Then for $\alpha$ small enough:
\begin{itemize}
\item $A_{\alpha}^{-1}(f^{-1})(\mathcal{O}^-) \Subset \mathcal{O}^-$,
\item $A_{\alpha}^{-1}(f^{-1}(\mathcal{O}^-_1)) \Subset \mathcal{O}^-_1$,  
\item $A_{\alpha}( f(\mathcal{O}^+)) \Subset \mathcal{O}^+$, 
\item $A_{\alpha}( f(\mathcal{O}^+_1)) \Subset \mathcal{O}^+_1$. 
\end{itemize}
Consider the element $g$ in $\mathrm{Orb}(f)$ defined as:
$$ g := A_{\alpha} \circ f \circ A_{\alpha}.$$
Observe that $I(g)= A_{\alpha}^{-1}(I(f))\subset A_{\alpha}^{-1}(H) \subset H$.  Hence we can assume 
that $I(g)\cap \overline{(\mathcal{O}^+ \cap \mathcal{O}^-)}=\varnothing$.  
 The following property are then satisfied:
\begin{itemize}
\item $g^{-1}(\mathcal{O}^-)  \Subset \mathcal{O}^-$ and $g(\mathcal{O}^+)\Subset \mathcal{O}^+$;
\item  $I(g)\cap \overline{(\mathcal{O}^+ \cap \mathcal{O}^-)}=\varnothing$.
\end{itemize}
Observe that the previous properties are stable under small perturbations. That defines a small open set $W_0$ in $W$ where the above conditions are satisfied and we can find a smooth probability measure $\mu$ with support in 
$W_0$ such that the above conditions are satisfied for $g_c=g\circ c$ with $c\in \mathrm{Supp}(\mu)$. 

Now, in order to prove Lemma \ref{pourk}, we choose for $\Omega_1$ any smooth form in $\C_1$ 
(for example the Fubini-Study form). 
As before, we take for $\Omega_s$ a smooth form in 
$\C_s$ with support in $\mathcal{O}^-$ and for $\beta_{k-s}$
a smooth form in $\C_{k-s}$ with support in 
$\mathcal{O}^+$.

In particular, for $c \in W_0$, $L^n_c(\Omega_s)\wedge  \Lambda^n_c(\beta_{k-s})$ is a probability measure with 
support in $\mathcal{O}^+\cap \mathcal{O}^-$. Lemma 2.3.5 in \cite{DS6} 
implies that there is a constant $C>0$ independent of $c$ such that
$$\|U_{ L_c(\Omega)}\|_{\C^1(\mathcal{O}^+_{s}\cap \mathcal{O}^-_{k-s})}\leq  C.$$
So arguing as above, we have that $k_n(c)$ is uniformly bounded for $c \in \mathrm{Supp}(\mu)$. That gives Lemma \ref{pourk} in
the case where $I$ is contained in a hyperplane.\\

We claim that in that example one also has that for any $x\in \P^k$,  $\log\mathrm{dist}(., x)$ is integrable with respect to $\nu_c$ for $c$ outside a pluripolar set (the distance being given by the Fubini-Study metric). The proof of that claim follows the lines of the previous one.  Choosing suitable coordinates, we can assume that $x=0 \in \Cc^k \subset \P^k$. Let $[z_0:\dots:z_{k-1}:t]$ denote the associated homogeneous coordinates on $\P^k$. We want to construct an example stable by perturbations for which  $\nu_{c}$ integrates $\log \mathrm{dist}( . , 0)$ with locally uniform estimates. 
The qpsh function $\log \|(z_0, \dots, z_{k-1})\| -\log \|(z_0, \dots, z_{k-1},t)\|$ is well defined so $dd^c \log \|(z_0, \dots, z_{k-1})\|$ is a well defined $(1,1)$ current in $\P^k$. Furthermore,  $\log \mathrm{dist}( . , 0) \in L^1(\nu_{c'})$ is equivalent to  $\log \|(z_0, \dots, z_{k-1})\| -\log \|(z_0, \dots, z_{k-1},t)\| \in L^1(\nu_{c'})$. 
Using super-potential theory (\cite[Lemma 4.2.8]{DS6}), it is enough to prove that 
$$\U_{T^+_{c', s}}(dd^c \log \|(z_0, \dots, z_{k-1})\| \wedge T^-_{c', k-s}) >L$$
for $c'$ in a small neighborhood of $c$. Observe that $dd^c \log \|(z_0, \dots, z_{k-1})\|$ is a well defined $(1,1)$ current in $\P^k$, since it is smooth outside a set of dimension $0$, its wedge product with any positive closed current is well defined (see \cite{dem}). We take $c\in W_0$ as in the previous example: 
\begin{itemize}
\item let $U_{c,s}$ denotes the Green quasi-potential of $L_c(\Omega_{s})$ of the previous section. Then $U_{c,s}$ is smooth (with locally uniform estimate near $c$) in $\mathcal{O}^-$.
\item  for all $n\geq 0$ and $c'$ near $c$, $\Lambda_{c'}^n(dd^c \log \|(z_0, \dots, z_{k-1})\wedge \Lambda^n_{c'}(\beta_{k-s}))$ is a well defined element of $\C_{k-s+1}$ with support in $\mathcal{O}^-$. It is $(f_{c'})_*$-admissible since $U_{c',s}$ is finite at that point. 
\item the functions
$$ k'_{m}(c'):=  \sum_n \U_{ L_{c'}(\Omega_s)}( \Lambda^n_{c'}( dd^c \log \|(z_0, \dots, z_{k-1})\| \wedge  \Lambda^m_{c'}(\beta_{k-s})))$$ 
satisfy $dd^c k_m =T^+_m-T^-_m$ with $\|T^\pm_m\| \leq C$ where $C$ does not depend on $m$.
\end{itemize}
So arguing as above, we deduce the claim. 

\begin{Remark} \rm
The parameters $c \in W_0$ give functions $f_c$ which are \emph{horizontal-like maps} in $\mathcal{O}^+ \cap \mathcal{O}^-$. Such maps were introduced by Dujardin in dimension 2 (see \cite{duj}) and have been extensively studied by Dinh, Nguyen and Sibony in  \cite{DS11, DinhNguyenSibony1}. In that last article, the authors prove in the inversible case that the measure $\nu_c$ is PB (of entropy $\log d_s$ and hyperbolic) that means that it integrates qpsh functions and in particular  $\log\mathrm{dist}(., x)$. 
\end{Remark}

\begin{Remark} \rm
We can extend the results of Theorem \ref{principal} to any map such that "there exists linear subspaces $E^+_s$ and $E^-_{k-s }$ of dimension $s$ and $k-s$ such that, up to a linear change of coordinates, the ball in $E^+_s$ of center $p=E^+_s\cap E^-_{k-s}$ and radius $\mathrm{dist}(p, f^{-1}(E_{k-s})\cap E^+_s)$ does not contain a point of $I(f)$".
Using that condition, we leave to the reader the proof of Theorem \ref{principal} for a map such that $\mathrm{dim}(I)=k-s$ and $\mathrm{Vol}(I)\leq k-s$.
 In general, that condition is not easy to verify and there is no reason for an arbitrary map to check it.
\end{Remark}

\noindent {\bf Question.}  The following question is natural in the settings of generic dynamics. Indeed, it is known to be false in the general case (see \cite{Dilgu}). For $c$ outside a pluripolar set, is the measure $\nu_c$ $PB$ (does it integrate $DSH$ functions)? If the answer was yes, one would deduce that the Lyapunov exponents are generically not $-\infty$ and that the measure $\nu_c$ does not charge pluripolar sets. \\ 

\noindent $\bullet$ {\bf Proof of point 5 in Theorem \ref{principal}}

Observe that for polynomials, one always have that $I$ is contained in the hyperplane at infinity. In the previous case, we have built an example using an element $A_\alpha$ that fixes $H$. When $H$ is the hyperplane at infinity, that means that $A_{\alpha}$ is an affine automorphism of $\Cc^k$. Since that example is stable under small perturbations in $\mathrm{Aff}(\Cc^k)$, we just have to compute the $dd^c$ of the different functions used in the previous part ($g_n$, $g'_{n,m}$, $k_n$). That is done exactly in the same way, observe that $W_1=\mathrm{Aff}(\Cc^k)$ is a Zariski dense open set in $\widetilde{W}_1\simeq \P^{k^2+k}$.\\

\noindent $\bullet$ {\bf Ergodicity and mixing}

Let $c$ be a generic parameter. It is natural to ask if the measure $\nu_c$ is mixing (or ergodic, but mixing is stronger). 
We are able to do so under an additional hypothesis : we need that $\nu_c$ does not charge $I'(f_c)$. The strategy is classical in complex analysis so we only sketch it:
\begin{enumerate}
\item one first show that the Green current $T^+_{s,c}$ is extremal in the sense that if $S\in \C_s$ is such that $S\leq T^+_{s,c}$ then $S=T^+_{s,c}$.
\item one proves that, for a smooth function $\varphi$, $\varphi \circ f_c^n T^+_{s,c}$ converges in the sense of currents to $c(\varphi) T^+_{s,c}$ where $c(\varphi)= \langle \varphi, \nu_c \rangle$ (at this point, one uses that the potentials of $T^+_{1,c}$ are integrable with respect to $\nu_c$).
\item one deduces that for $\psi$ smooth, we have $\lim_n \langle \varphi \circ f^n \psi , \nu_c \rangle= \langle \varphi  , \nu_c \rangle \langle \psi  , \nu_c \rangle$. This would be true by the above if $T^-_{s,c}$ was smooth and one proceed by approximations (we need here that $\nu_c$ does not charge $I'(f_c)$). The mixing is proved. 
\end{enumerate}

\noindent $\bullet$ {\bf Hyperbolicity of the homogeneous extension and hyperbolicity of the map}

Assume now that $f$ is a dominating meromorphic map of $\P^k$. We can write it in homogeneous coordinates as
$f=[P_0:\dots:P_k]$ where the $P_i$ are relatively prime homogeneous polynomials of degree $d$. We consider the polynomial map of $\Cc^{k+1}$ defined as:
$$\widetilde{f} =  (P_0,\dots,P_k).$$ 
Its extension to $\P^{k+1}$ (still denoted as $\widetilde{f}$) has its indeterminacy set contained in $H$, the hyperplane at infinity. Hence, it satisfies the above conditions.  Let $[z_0:\dots:z_k:t]$ be the homogeneous coordinates on $\P^{k+1}$. Let $(\widetilde{d_i})_{i=0..k+1}$ be the generic dynamical degree in the orbit of $\widetilde{f}$. An easy computation gives:
$$ \widetilde{d_0}=1, \ \widetilde{d_i}= d\times d_{i-1} \ \mathrm{for} \ i\neq 0. $$
In particular, we can apply point 5 of Theorem \ref{principal} to $\widetilde{f}$. \\

Assume furthermore that $\widetilde{f}$ is in fact a hyperbolic map in the sense that it satisfies Theorem \ref{principal}. In other words, the parameter $\mathrm{Id} \in \mathrm{Aut}(\P^{k+1})$ is not in the pluripolar set where we cannot apply the Theorem.  Let $\widetilde{\nu}$ denote the measure of maximal entropy constructed for $\widetilde{f}$. Observe that since $0$ is an attractive fixed point, it does not belong to the support of $\widetilde{\nu}$, hence $\log \mathrm{dist}(x, 0) \in L^1(\widetilde{\nu})$. \\

The mapping $\widetilde{f}=[f: t^d]$ is a skew-product over $f$: if $\pi$ denotes the (meromorphic) projection from $\P^{k+1}$ to $\P^k$ defined by $\pi ([z_0:\dots:z_k:t]) =[z_0:\dots:z_k]$ then $f \circ \pi = \pi\circ \widetilde{f}$. In fact, as $\widetilde{\nu}$ does not charge $0$ (since it integrates $\log \mathrm{dist}(., 0)$), we can work instead in the birational model $\P^k \times \P^1$ where the map $\pi$ is holomorphic. Let $\nu':= \pi_* \widetilde{\nu}$. We claim that:
\begin{theorem}
Assume that $\widetilde{f}$ is as above, then the measure $\nu'$ is a hyperbolic measure of maximal entropy
$\log{d_s}$. Assume that the measure $\nu= T_s^+ \wedge T_s^-$ is well defined, then $\nu$ is also a hyperbolic measure of maximal entropy
$\log{d_s}$
\end{theorem}
\noindent \emph{Proof.} Using Proposition 3.5 in \cite{L.W} gives
	$$h_{\widetilde{\nu}}(\widetilde{f}) \leq h_{\nu'}(f) + \int_{\P^k=H} h(\widetilde{f}, \pi^{-1}(y)) d\nu'(y),$$
where $h(\widetilde{f}, \pi^{-1}(y))$ is the topological entropy of $\widetilde{f}$ relative to the set $\pi^{-1}(y)$. Observe that in \cite{L.W}, the mappings $\widetilde{f}$ and $f$ are continuous but that hypothesis is not needed for that inequality. On the other hand, on $\pi^{-1}(y)\simeq \P^1$ the mappings $\widetilde{f}_{y}:=(\widetilde{f})_{|\pi^{-1}(y)}$ are holomorphic maps of degree either $d$ or $0$ (that happens when $y\in I(f)$). Then  $h(\widetilde{f}, \pi^{-1}(y))$ is the entropy of the sequence $(\widetilde{f}_{y_n})_n$ where $y_n=f^n(y)$ (see \cite{K.L} for definitions). Gromov's arguments on $lov$ (see \cite{Gro}) still apply in that setting and one gets that  $h(\widetilde{f}, \pi^{-1}(y))\leq \log d$. In particular, we deduce:
$$\log d_s + \log d\leq h_{\nu'}(f)+ \log d.$$
In other words, $h_{\nu'}(f)\geq \log {
d_s}$. As the other inequality always stands  (see \cite{DS9}), that gives $h_{\nu'}(f)= \log {d_s}$.  \\

 Since $ \log  \|(z_0,\dots,z_k)\|- \log  \|(z_0,\dots,z_k,t)\| \in L^1(\widetilde{\nu})$, invariance of $\widetilde{\nu}$ implies that
\begin{align*}  
&\log  \|f(z_0,\dots,z_k)\|-d\log  \|(z_0,\dots,z_k)\|  +  \\
&d\log  \|(z_0,\dots,z_k )\| -  d\log  \|(z_0,\dots,z_k,t)\|+ \\
&  d\log  \|(z_0,\dots,z_k,t)\| - \log  \|\widetilde{f}(z_0,\dots,z_k,t)\| \in L^1(\widetilde{\nu}). 
\end{align*}
We have that $d\log  \|(z_0,\dots,z_k,t)\| - \log  \|\widetilde{f}(z_0,\dots,z_k,t)\|\in L^1(\widetilde{\nu})$ by hypothesis and $d\log  \|(z_0,\dots,z_k )\| -  d\log  \|(z_0,\dots,z_k,t)\| \in L^1(\widetilde{\nu}) $, thus: 
\begin{align*}
\log  \|f(z_0,\dots,z_k)\|-d\log  \|(z_0,\dots,z_k)\| = \\
\pi^*(\log  \|f(z_0,\dots,z_k)\|-d\log  \|(z_0,\dots,z_k)\|) \in L^1(\widetilde{\nu}).
\end{align*}
We deduce that $\nu'$ integrates  a quasi-potential of $f^*(\omega)$ hence $\log \mathrm{dist}(.,I)$. De Thélin's Theorem can be applied and we deduce the hyperbolicity of $\nu'$. \\

 By continuity of $\pi_*$, one has that 
$$\nu'= \lim_{n\to \infty} \pi_*( \frac{1}{(d \times d_s)^n}(\widetilde{f}^n)^*(\Omega_{s+1}) \wedge  \frac{1}{(d \times d_s)^n}(\widetilde{f}^n)_*(\Omega_{k-s}))$$
 where $\Omega_{s+1}$ and $\Omega_{k-s}$ are smooth elements of $\C_{s+1}(\P^{k+1})$ and $\C_{k-s}(\P^{k+1})$. In particular, we choose $\Omega_{s+1}= \pi^*(\omega^s) \wedge \Omega_1$ where $\omega$ is the Fubini-Study form on $\P^k$ and $\Omega_1$ is a smooth $(1,1)$ form with support disjoint from $0$ (observe that $\pi^*(\omega^s)$ is smooth away from $0$). Then we have that  $(\widetilde{f}^n)^*(\Omega_{s+1})= \pi^*(f^n)^*(\omega^s)\wedge  (\widetilde{f}^n)^*(\Omega_1)$. Thus:
\begin{align*} 
&\pi_*( \frac{1}{(d\times d_s)^n}(\widetilde{f}^n)^*(\Omega_{s+1}) \wedge  \frac{1}{(d\times d_s)^n}(\widetilde{f}^n)_*(\Omega_{k-s}))= \\ 
&\frac{1}{d_s^n}(f^n)^*(\omega^{s}) \wedge \pi_*( \frac{1}{d^n}(\widetilde{f}^n)^*(\Omega_1) \wedge  \frac{1}{(d\times d_s)^n}(\widetilde{f}^n)_*(\Omega_{k-s}) ).
\end{align*}
Now, $\frac{1}{d_s^n}(f^n)^*(\omega^s)$ converges in the Hartogs' sense to the Green current $T^+_{s}$ of $f$. Let $\widetilde{T}^+_{1}$ and $\widetilde{T}^-_{k-s}$ be the Green currents of $\widetilde{f}$. They are wedgeable by hypothesis.  In particular, 
$\pi_*( \frac{1}{d^n}(\widetilde{f}^n)^*(\Omega_1) \wedge  \frac{1}{(d \times d_s)^n}(\widetilde{f}^n)_*(\Omega_{k-s}) )$ converges in the Hartogs' sense to $\pi_*( \widetilde{T}^+_{1} \wedge \widetilde{T}^-_{k-s})$. One easily checks that it defines an $f_*$-invariant current in $\C_{k-s}(\P^k)$.
 As $T^-_{s}$ is the more H-regular invariant current, we deduce that it is more H-regular than  $\pi_*( \widetilde{T}^+_{1} \wedge \widetilde{T}^-_{k-s})$. Thus $\nu$ is more H-regular than $\nu'$ and it particular, $\nu$ integrates $\log \mathrm{dist}(., I)$. We can then apply as above  Theorem \ref{ENTROPY} and Corollary 3 in \cite{DT1} to compute the entropy and prove the hyperbolicity of $\nu$. \hfill $\Box$ \hfill   \\

\noindent {\bf Question.}  Unicity of the measure of maximal entropy is expected so it would nice to prove that $\nu'=\nu$. \\

\noindent Gabriel Vigny, LAMFA - UMR 7352, \\ 
U. P. J. V. 33, rue Saint-Leu, 80039 Amiens, France. \\
\noindent Email: gabriel.vigny@u-picardie.fr

\end{document}